\documentclass[aip,graphicx, preprint]{revtex4-1}

\pdfoutput=1

\usepackage{amsfonts,amssymb,amsmath,amsthm}
\usepackage{fullpage}
\usepackage{bm}
\usepackage{bbm}
\usepackage{enumerate}
\usepackage{dsfont}
\usepackage{pxfonts}
\usepackage{xcolor}
\usepackage[colorlinks,
pdffitwindow=false,
plainpages=false,
pdfpagelabels=true,
pdfpagemode=UseOutlines,
pdfpagelayout=SinglePage,
bookmarks=false,
colorlinks=true,
hyperfootnotes=false,
linkcolor=blue,
citecolor=green!50!black]{hyperref}

\usepackage{tikz}
\usetikzlibrary{decorations.markings} 
\usetikzlibrary{arrows,decorations.pathmorphing,backgrounds,fit,positioning,shapes.symbols,chains}
\usepackage{pgfplots}
\usepgfplotslibrary{fillbetween}

\newcommand{\R}{\mathbb{R}}

\newcommand{\set}[1]{\mathbb{#1}}		

\newcommand{\ep}{\varepsilon}
\newcommand{\Cl}{\mathrm{cl}} 

\graphicspath{{./graphics/}{./graphics/double_gyre/}{./graphics/1dSys/}{./graphics/aviso/}}

\newtheoremstyle{custom}{3pt}{3pt}{}{}{\bfseries}{:}{.5em}{}
\theoremstyle{custom}
\newtheorem{example}    {Example}

\newtheorem{remark} 		[example]{Remark}
\newtheorem{proposition}[example]{Proposition}
\newtheorem{assumption}[example]{Assumption}

\definecolor{FUblue}{RGB}{0,51,102}
\definecolor{FUgreen}{RGB}{153,204,0}
\definecolor{FUred}{RGB}{204,0,0}

\draft 

\begin{document}


\title{Network Measures of Mixing} 



\author{Ralf Banisch}
\email[]{ralf.banisch@fu-berlin.de}
\affiliation{Institute of Mathematics, Freie Universit\"at Berlin, 14195 Berlin, Germany. }
\author{P\'eter Koltai}
\email[]{peter.koltai{@}fu-berlin.de}
\affiliation{Institute of Mathematics, Freie Universit\"at Berlin, 14195 Berlin, Germany. }
\author{Kathrin Padberg-Gehle}
\email[]{padberg@leuphana.de}
\thanks{}
\affiliation{Institute of Mathematics and its Didactics, Leuphana Universit\"at L\"uneburg, Universit\"atsallee 1, 21335 L\"uneburg, Germany. }


\date{\today}

\begin{abstract}
Transport and mixing processes in fluid flows can be studied directly from Lagrangian trajectory data, such as obtained from particle tracking experiments. Recent work in this context highlights the application of graph-based approaches, where trajectories serve as nodes and some similarity or distance measure between them is employed to build a (possibly weighted) network, which is then analyzed using spectral methods.  Here, we consider the simplest case of an unweighted, undirected network and analytically relate local network measures such as node degree or clustering coefficient to flow structures. In particular, we use these local measures to divide the family of trajectories into groups of similar dynamical behavior via manifold learning methods.
\end{abstract}


\maketitle 


\begin{quotation}
Coherent structures play a ubiquitous role in mass transport and mixing in time-dependent dynamical systems. While most of the established identification methods require the knowledge of the underlying dynamical system or at least high-resolution trajectory data, this information may not be available in real-world scenarios. Network-based approaches applied to Lagrangian trajectories have been shown to be still successful in the analysis of transport and mixing processes even when data are sparse and incomplete. We focus on the local network measures of such trajectory networks and relate them to flow structures of the underlying system.
Using manifold learning algorithms allows us to classify the phase space into regions of different dynamical behavior, complementing the frequently used spectral studies of flow networks. 
\end{quotation}

\section{Introduction}
To set the scene for our approach, suppose we are given a nonautonomous ordinary differential equation (ODE)
\begin{equation}\label{eq:dgl}
\dot{x}(t)=f(t, x(t))
\end{equation}
with state $x\in  \set{M} \subset \mathbb{R}^d$, $d\in \mathbb{N}$, time $t \in \mathbb{R}$ and sufficiently smooth right-hand side 
$f$ to ensure existence and uniqueness of solutions $x(t):=x(t; t_0, x_0)$ to initial values $(t_0,x_0)$. We may interpret $f$ as a 
velocity field of a fluid flow and $x(t)$ as the position of an ideal particle in that flow at time $t$. 
In this case, $\set{M} \subset \mathbb{R}^2$ or  $\mathbb{R}^3$. The time-parameterized family $(x(t))_t$ is referred to as a Lagrangian particle trajectory.
We are interested in detecting coherent flow structures from the given trajectory data, that is, time-dependent regions in $\set{M}$ that either inhibit or enhance 
transport and mixing processes of the underlying flow. 
In our setting, transport refers to the bulk movement of particles between different regions of the domain. Coherent structures are characterized by trapping particles for long times and in that way they determine the transport pathways of the underlying flow. Mixing refers to advective mixing (or stirring) by stretching and folding of fluid parcels and compared to transport it describes a local flow property. Coherent structures resist mixing with neighboring regions of the flow and thus they are characterized by having other mixing properties than these surrounding regions. 

The mathematical definition and numerical study of coherent structures has been an area of intense research over the last two decades. In particular, different probabilistic and geometric methods have been proposed. These approaches have been discussed and systematically compared in a number of studies\cite{FrPa09,Allshouse_Peacock_2015, HaFa_etal17}. Geometric concepts aim at defining the boundaries between coherent sets, i.e., codimension-1 material surfaces in the flow that can be characterized by variational criteria\cite{Haller_Rev_2015}. Probabilistic methods are tailored to identify sets that are minimally dispersive while moving with the flow. Here the main theoretical tools are transfer operators, i.e., linear Markov operators that describe the motion of probability densities under the action of the nonlinear, time-dependent flow\cite{FPG14}. 
Both the geometric and the probabilistic approach require high resolution trajectory data from \eqref{eq:dgl}, that is, from a dense grid of initial conditions. This can be prohibitively expensive in complex systems, such as turbulent flows. Moreover, when the particle trajectories are obtained directly from measurements (e.g., from particle tracking experiments), then the Lagrangian data under consideration may even be sparse and incomplete.

To overcome these problems, different computational methods have been proposed to identify coherent behavior in flows directly from Lagrangian trajectory data. One of the earliest attempts is the braiding approach\cite{AT12}, where trajectories are classified according to their  intertwining pattern in space-time. This method is mathematically sound, but computationally demanding and currently restricted to two-dimensional flows. Other trajectory-based approaches use time-integrated quantities along trajectories\cite{levnajic2010ergodic,MeLoFoHo10,BM12,Mancho2013}.  Finally, there are attempts to reconstruct the transfer operator from limited amount of trajectory data\cite{Williams_et_al_2015} as well as the dynamic Laplacian\cite{FrJu18}, which was recently introduced to study coherent sets as sets that keep an optimal boundary to volume ratio when evolved by the dynamics\cite{Froyland_2015}.

Recent works focus on the use of spatio-temporal clustering algorithms\cite{Froyland_Padberg_2015,Hadjighasem_et_al_2016,Banisch_Koltai_2017,Schlueter_Dabiri_2016,PGSc17}. There, the aim is to identify coherent sets as groups of trajectories that remain close and/or behave similarly in the time span under investigation. All these methods can deal with sparse and incomplete trajectory data and their applicability has been demonstrated in several example systems.  

Here, we revisit the framework introduced in Ref. \onlinecite{PGSc17}, which is based on an unweighted, undirected network with the trajectories serving as nodes.  A link is established between two nodes if the respective trajectories come close to each other at least once in the course of time. This construction is similar in spirit to the concept of recurrence networks\cite{Donner_et_al_2010a,Donner_et_al_2010b}, but here in a spatio-temporal setting. 

We note that the discretized transfer operator has also been viewed and treated as a network\cite{dellnitz_preis_03,dellnitz_etal_05,Padberg_et_al_2009,Lindner_Donner_2017,Ser-Giacomi_et_al_2015, RoSeHe17}. A recent review\cite{Donner2019} addresses the different constructions of flow networks and their analyses.

In previous work\cite{PGSc17}, we have introduced the construction and have mainly considered spectral properties of our trajectory-based undirected and unweighted flow network, which allowed us to compare our approach with related spectral concepts\cite{Hadjighasem_et_al_2016,Banisch_Koltai_2017}. We note that in these works the search for coherent sets is termed as a community detection problem for the resulting network and solved by a normalized cut method\cite{shimalik} . 

In the present paper, our focus will be on the application and interpretation of \emph{local network measures} such as node degrees or clustering coefficients.  
These and other quantities have been considered in previous work on recurrence networks\cite{Donner_et_al_2010a}, where the authors could link the network measures to properties of the underlying dynamical system.  In weighted, directed networks obtained from discretized transfer operators the in- and out-degrees where found to highlight hyperbolic regions in the underlying flow\cite{Ser-Giacomi_et_al_2015,Lindner_Donner_2017}, whereas maxima of the local clustering coefficient could be related to regular or periodic behavior\cite{RoSeHe17}.
Similar properties appear to hold for the trajectory-based undirected and unweighted flow network, as demonstrated in example systems\cite{PGSc17}.

Here, for the first time, we will draw an explicit analytical connection between these network measures and underlying flow structures. In particular, we will estimate the node degree in terms of the finite-time Lyapunov exponent, an established quantity to measure stretching, and give some geometrical interpretation of the local clustering coefficient. 
Moreover, we will carry out an empirical flow classification based on further network measures using manifold learning methods. These complement the spectral approaches\cite{Hadjighasem_et_al_2016,Banisch_Koltai_2017,Schlueter_Dabiri_2016,PGSc17}.

The paper is organized as follows: In section \ref{sec:network} we review the construction of the simple trajectory-based network\cite{PGSc17} as well as standard network local measures. In section \ref{sec:degree} we will establish an analytical connection of some local network measures to the corresponding phase space structures. In particular, we will give analytical estimates of the node degree and of the local clustering coefficient. Further network measures are discussed in section \ref{sec:furthermeas}. In section \ref{sec:numerics} we will numerically demonstrate the estimates of section \ref{sec:degree} as well as an empirical network-based flow classification in a number of example systems, including the double-gyre flow and a real ocean surface flow from the AVISO data set.

\section{A trajectory-based network}\label{sec:network}
\subsection{Construction of the network} \label{ssec:network_constr}

Suppose we are given $N$ trajectories from a flow simulation (i.e., as numerical solutions to \eqref{eq:dgl}) or from a particle tracking experiment. In practice, the particle positions may be given at discrete times $\{0, 1,\ldots, T\}$. We denote the trajectories by $x_{i}$, $i=1, \ldots, N$, 
and the positions at a certain time instance $t=0,\ldots,T$ by $x_{i,t}\in \mathbb{R}^d$. We now fix some $\ep>0$.
and set up a network with $x_1, \ldots, x_N$ as nodes. We link two nodes $x_i$ and $x_j$ if the respective trajectories come $\ep$-close to each other at least once in the course of time~$t\le T$. Then, a symmetric adjacency matrix $A\in \{0,1\}^{N \times N}$ describes the network, with
\begin{equation}
A_{ij}=
\begin{cases} \max\limits_{0 \leq t \leq T} \mathbf{1}_{\set{B}_{\ep}(x_{i,t})}(x_{j,t}),& i\neq j \\ 
0, & i=j \end{cases}, \label{eq:A}
\end{equation} 
where $\mathbf{1}_{\set{B}}$ denotes the indicator function of a set $\set{B} \subset \mathbb{R}^d$.
So $A_{ij}=1$, that is, a link is established between trajectories $x_i$ and $x_j$, if and only if at one or more time instances $t$, $x_{j,t}$ can be found in an $\ep$-ball $\set{B}_{\ep}(x_{i,t})$ centered at $x_{i,t}$ and thus the trajectories $x_i$ and $x_j$ have come $\ep$-close.

Naturally, other constructions are possible as well, e.g., setting $A_{ij}$ to the number of $\ep$-close encounters and similar. We will restrict our attention to the current one, as we consider it to be the simplest in terms of carrying the least quantitative information.

By an appropriate choice of $\ep$ one ensures that the network defined by \eqref{eq:A} is connected and we will only consider connected networks from now on. Of course, the network topology depends also on other parameters, such as the particle density, the time-resolution and the length of the trajectories. However, as these parameters are determined by the data, the only free parameter is $\ep$, which can be reasonably chosen in relation to the particle density\cite{Donner_et_al_2010b,PGSc17}.

\subsection{Network analysis}  \label{ssec:network_ana}
The resulting network can be studied globally using the adjacency matrix. In particular, a normalized cut problem\cite{shimalik} is solved by considering leading eigenvectors of the generalized eigenvalue problem $Lv= \lambda  Dv$, where $D$ denotes the degree matrix (diagonal matrix with node degrees on the diagonal, i.e., $D_{ii}=\sum_{j}A_{ij}$, $i=1, \ldots, N$). This spectral approach leads to the identification of clusters in the network that correspond to coherent sets of the underlying flow\cite{Hadjighasem_et_al_2016,Banisch_Koltai_2017,Schlueter_Dabiri_2016,PGSc17}, that is, mobile regions in $M$ that do not freely mix with the surrounding phase space regions. 

As an alternative to these global, spectral approaches, we will attempt to tie in \emph{standard local network measures}\cite{newman2003} with quantitative and qualitative dynamical behavior of the system. So far such relations have been studied mainly heuristically\cite{Padberg_et_al_2009,Donner_et_al_2010a,Ser-Giacomi_et_al_2015,Lindner_Donner_2017}, but here we will draw some explicit analytical connections.   

In the following, we will review some of the most frequently considered network measures.

The \emph{node-degree}
\begin{equation}
d_i =\sum_j A_{ij}, \;\;\ i=1, \ldots, N,\label{eq:k}
\end{equation}
counts how many links are connected to node~$x_i$. Similarly, 
\begin{equation}
\langle d \rangle_{nn, i} = \frac{\sum_j A_{ij}d_j}{d_i}, \;\; i=1, \ldots, N, \label{eq:knn}
\end{equation}
denotes the \emph{average node degree} of the neighbors of a node $x_i$.  

The \emph{local clustering coefficient} of a node $x_i$ quantifies how tightly connected the subgraph induced by this node and its neighbors is (that is, how close this subgraph is to a complete graph). It is defined by
\begin{equation}
c_i= \frac{\# \mbox{ triangles connected to } x_i} { \# \mbox{ triples centered around } x_i} =\frac{(A^3)_{ii}}{d_i(d_i-1)}, \label{eq:cc}
\end{equation}
$ i=1,\ldots, N$. 

Many further network measures exist in the literature, we will only mention a few here:

The (normalized) \emph{closeness centrality} (or \emph{closeness}) of a node in a network is given by its reciprocal mean distance to all other nodes in the network, i.e.,
\[
\Cl_i = \frac{N}{\sum_j \mathrm{dist}(x_i,x_j)},
\]
where $\mathrm{dist}(x_i,x_j)$ denotes the distance (or shortest path length) between nodes $x_i$ and~$x_j$ in the graph. 

\emph{Betweenness centrality} is related to that quantity and measures how many shortest paths in the network traverse a certain node. It is given by
\[
b_i = \sum_{j,k \neq i}^N \frac{p_{jk}(x_i)}{p_{jk}},
\]

where $p_{jk}(x_i)$ denotes the number of shortest paths between nodes $x_j$ and $x_k$ that include $x_i$ and $p_{jk}$ the total number of shortest paths between $x_j$ and $x_k$.

\section{Analytical estimates of local network measures}
\label{sec:degree}
In a number of previous studies in recurrence networks\cite{Donner_et_al_2010a}, transfer operator-based\cite{Padberg_et_al_2009,Ser-Giacomi_et_al_2015,Lindner_Donner_2017} as well as trajectory-based networks\cite{PGSc17} it has been observed that high values of the node degree can be related to regions of strong stretching, indicating hyperbolic behavior of the underlying flow. 
Moreover, high values of the local clustering coefficient have been related to regular/elliptic dynamics and periodic behavior\cite{Donner_et_al_2010a,RoSeHe17}.  Here we will establish an analytic connection between these local network measures and the finite-time Lyapunov exponent, a frequently used indicator of hyperbolic dynamics. 

Although we consider time-continuous dynamics for our analytic investigations, the results here \emph{qualitatively} carry over to discrete-time systems as well, by replacing the corresponding discrete-time analogs of the objects in consideration. The quantitative result, Proposition~\ref{prop:vol_galaxy} below, can be retained from the discrete-time system in the limit of vanishing sampling step size, if the discrete-time system is the finite-time flow map of a continuous-time system.


\subsection{Preliminaries}

Given a time-dependent flow~$\phi(s,t)$ generated by the ODE \eqref{eq:dgl},
i.e., $x(t) = \phi(s,t)x(s)$ that solves $\dot{x}(t) = f(t,x(t))$, let us assume that this flow generates our trajectory data. Herein,~$f:\R\times\R^d\to\R^d$ is a sufficiently smooth vector field.
For simplicity, we assume that the flow is volume-preserving for all times, i.e.,~$\text{div}_x\,f(t,\cdot) \equiv 0$ for all~$t$, where~$\text{div}_x$ is the divergence operator for functions mapping~$\R^d$ to itself.

Let~$W(s,t)\in\R^{d\times d}$ be the \emph{fundamental matrix} (Wronskian matrix), defined by
\begin{equation}
W(s,t) = D_x \left(\phi(s,t)x\right)\,,
\label{eq:Wronski}
\end{equation}
where~$D_x$ denotes the derivative with respect to~$x$. Hence, the fundamental matrix is the derivative of the flow with respect to its initial condition. It is a first-order (linear) approximation of perturbation propagation along trajectories, since it holds that if~$\tilde{x}$ satisfies~\eqref{eq:dgl} with~$\tilde{x}(s) = x(s)+\delta x$, then
\begin{equation}
\tilde{x}(t) = x(t) + W(s,t)\delta x + \mathcal{O}(\|\delta x\|^2)\,,
\label{eq:perturbation}
\end{equation}
as $\|\delta x\|\to 0$. The fundamental matrix satisfies the initial value problem\footnote{Sometimes \emph{any} solution $\tilde{W}(s,t)$ of~\eqref{eq:vareq} for an arbitrary initial condition is called fundamental matrix. In this case, our definition from~\eqref{eq:Wronski} is recovered by~$W(s,t) = \tilde{W}(s,t)\tilde{W}(s,s)^{-1}$.}
\begin{equation}
\dot{W}(s,t) = D_xf(t,x(t))\, W(s,t),\qquad W(s,s) = I\,,
\label{eq:vareq}
\end{equation}
where the time-derivative on the left-hand side is with respect to~$t$, and~$I\in\R^{d\times d}$ denotes the identity matrix.
We note the following properties:
\begin{proposition} \label{prop:det1}
As the vector field~$f$ is divergence-free, we have $\det (W(s,t))=1$ for every $x$ and $s\le t$.
\footnote{This follows from the chain rule, from the fact that $\frac{d}{dt}\det (A(t)) = \text{tr}(\frac{1}{\det (A(t))}A(t)^{-1} \frac{d}{dt}A(t))$ and that the trace is invariant under similarity transformations.}
\end{proposition}
The finite-time Lyapunov exponent (FTLE) for $x=x(s)$ is defined as
\begin{equation}\label{eq:ftle}
\Lambda(x, s; t)= \frac{1}{|t-s|} \log\left( \sigma_1(W(s,t))\right),
\end{equation}
 where $\sigma_1 (W(s,t))$ denotes the largest singular value of the matrix $W(s,t)$.

Let us now turn our attention to the properties of the dynamic neighborhood network defined by~\eqref{eq:A}.
Fix some~$\ep>0$.
Let us recall our initial setting again, where we have some sampling points and their trajectories in a volume-preserving flow. Assuming the initial distribution of the sampling points being uniform (in some spatial domain of interest), and there are sufficiently many of them, the number of sampling points in the $\ep$-neighborhood of any point is proportional to the volume of the $\ep$-ball around that point.\footnote{This statement is made rigorous by considering independently uniformly sampled points, then their relative ratio in the $\ep$-neighborhood converges to the relative volume of this neighborhood almost surely as the number of samples goes to infinity. This is the consequence of the law of large numbers, exploited in Monte Carlo methods.} 
Denote 
\begin{equation} \label{eq:Geps}
\begin{aligned}
\set{G}_{\ep} (x)  & = \bigcup_{t=t_0,\ldots,t_{T-1}} \left\{\tilde{x}\,\big\vert\, \phi(t_0,t)\tilde{x} \in \set{B}_{\ep}(x(t))\right\} \\
&=  \bigcup_{t=t_0,\ldots,t_{T-1}}  \phi(t_0,t)^{-1} \set{B}_{\ep}(x(t)),
\end{aligned}
\end{equation}
where~$x(t)$ solves~\eqref{eq:dgl} with~$x(t_0)=x_{t_0}^i$ and $\phi(t_0,t)^{-1}:=\phi(t, t_0)$ denotes the time-reversed flow from time $t$ back to $t_0$. Thus, $\set{G}_{\ep}(x)$ denotes the set of all initial states that are $\ep$-close to the trajectory starting at $x$ at some $t_i$ for $i = 0, 1, \ldots, T$.

\subsection{The degree.}
\label{ssec:degree}
As the initial distribution is uniform, in the case of many data points (more precisely, in the $N\to\infty$ limit) the Lebesgue volume of a set is proportional to the number of data points lying in it.
Thus, the degree $d_i$ of node~$x_i$
in the network given by the adjacency matrix~\eqref{eq:A} (the row sum of $A$) is proportional to the set of initial conditions that get $\ep$-close to~$x_i(t)$ at some time in the observation interval, i.e.,
\begin{equation} \label{eq:degreepropvolume}
d_i=\sum_{j=1}^N A_{ij} \propto \text{vol}(\set{G}_{\ep} (x_i)).
  \end{equation}
Here, the proportionality holds for large $N$ and with the usual Monte Carlo error of~$\mathcal{O}(N^{-1/2})$.
We will now estimate the volume in~\eqref{eq:degreepropvolume} to first order in~$\ep$ using the linear perturbation propagation relation~\eqref{eq:perturbation}.
For a fixed~$t$ let~$\tilde{x}$ be such that~$y=\phi(t_0,t)\tilde{x}\in B_\ep(x_i(t))$, and~$\delta y = y - x_i(t)$. Then, by~\eqref{eq:perturbation} we obtain
\begin{equation}
\|\tilde{x}-x_i(t_0)\| = \|W(t_0,t)^{-1}\delta y\| + \mathcal{O}(\ep^2)\,.
\label{eq:pullbackdist}
\end{equation}
Note that~$W(s,t)$ is always invertible due to Proposition~\ref{prop:det1}.

We have from~\eqref{eq:pullbackdist} that
\begin{equation}
d_i \,\dot{\propto}\ \text{vol}\Big(x_i(t_0) + \bigcup_{t=t_0,\ldots,t_{T-1}} W(t_0,t)^{-1} \set{B}_{\ep}(0)\Big)\,,
\label{eq:pullbackballs}
\end{equation}
where~$\dot{\propto}$ denotes proportionality up to errors of size~$\mathcal{O}(\ep^2)$. Note that the translation by~$x_i(t_0)$ does not change the volume and could be omitted. The sets~$W(t_0,t)^{-1} \set{B}_{\ep}(0)$ are ellipsoids with semi-axes of length~$\sigma_1^{-1},\ldots,\sigma_d^{-1}$, the reciprocal singular values of the matrix~$W(t_0,t)$, and these semi-axes are aligned with the corresponding right singular vectors.

We now restrict our considerations to the case of a two-dimensional area-preserving flow with states $(x,y) \in \R^2$.

\paragraph{Time-invariant singular vectors.}

For~$d=2$, we have by Proposition~\ref{prop:det1} that~$\sigma_1\sigma_2=1$. We will make the following simplifying assumption.
\begin{assumption} \label{ass:alignment}
Let~$\sigma(t)$ denote the larger of the both singular values of~$W(t_0,t)$, and $v(t)\in\R^2$ the corresponding right singular vector. We assume that~$v(t)=v$ is independent of~$t$.
\end{assumption}
This means that the other right singular vector is also independent of~$t$. In other words, we assume, that the direction at initial time which undergoes the largest stretching (and also the one undergoing the largest squeezing) is independent of~$t$.

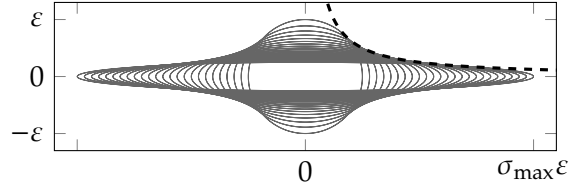
\begin{figure}[h]
\centering

\centering 
\begin{tikzpicture}[scale = 1]

\begin{axis}[
	xtick = {-4,0,4},
    xticklabels = {{},0,$\sigma_{\text{max}}\varepsilon$ },
    ytick = {-1,0,1},
    yticklabels = {$-\varepsilon$,0,$\varepsilon$},
    xmin = -4, xmax = 4,
	ymin = -1, ymax = 1,
	enlarge x limits=0.05,
	enlarge y limits=0.15,
	trig format plots=rad,
	variable=t,
	samples = 200,
	width = 0.5\textwidth,
	axis equal image
]

\foreach \r in {8,...,32}
\addplot[domain=-2*pi:2*pi, thin, black!60]
({\r/8*cos(t)}, {8/\r*sin(t)});

\addplot[domain=0:4.5, dashed, very thick]
{1/(2*t)};

\end{axis}
\end{tikzpicture}

\caption{Linearized pullbacks of $\ep$-circles. Dashed line: the graph of $y = \frac{1}{2x}$, touching the union of the pullback circles for $x\in[\frac{1}{\sqrt{2}},\frac{\sigma_{\text{max}}}{\sqrt{2}}]$. In this visualization $\sigma_{\text{max}} = 4$.
}
\label{fig:ellipseunion}
\end{figure}

In order to simplify notation, we fix initial time $t_0=0$ and final time~$T>0$. Let us consider the set
\[
\set{S}(0,T) = \bigcup_{t\in [0,T]} W(0,t)^{-1} \set{B}_{\ep}(0)\,,
\]
i.e., a union over continuous time. We would like to compute the two-dimensional volume of~$\set{S}(0,T)$. It is a union of ellipses with major semi-axes of length~$\sigma\in [0,\sigma_{\text{max}}]$, where we also assume~$\sigma_{\text{max}}=\sigma(T)$, i.e., the largest stretching appears at final time. Without loss, the larger semi-axis is assumed to be aligned with the $x$-axis. The resulting set~$\set{S}(0,T)$ is depicted in Figure~\ref{fig:ellipseunion}. Next we will derive an analytic formula for the volume of $\set{S}(0,T)$. Those readers not interested in the details of the derivation can skip to the result in Proposition~\ref{prop:vol_galaxy}.

\begin{figure*}[htb]
\centering

\begin{tikzpicture}[scale = 1]
    \begin{axis}[
    enlarge x limits=0.05,
	enlarge y limits=0.15,
    xtick       = {0,0.7071, 2.8284,4},
    xticklabels = {$0$, $\frac{\ep}{\sqrt{2}}$, $\frac{\sigma(T)\ep}{\sqrt{2}}$, $\sigma(T)\ep$ },
    ytick = \empty,
    samples = 200,
    xmin = 0, xmax = 4,
	ymin = 0, ymax = 1,
	width = 0.85\textwidth,
	axis equal image
    ]

    \addplot[name path=zero, domain= 0:4, very thin, black] {0};

	 \addplot[name path=circ, domain= 0:1, thick, dashed, gray] {sqrt(1-\x^2)};
    \addplot[name path=circ, domain= 0:1/sqrt(2), thick, black] {sqrt(1-\x^2)};
    \addplot[
        fill=FUred,
        fill opacity=0.2
    ]
    fill between[
        of=zero and circ,
        soft clip={domain=0:0.7071},
    ];
    \addplot[name path=hyp, domain= 1/sqrt(2):4/sqrt(2), thick, black] {1/(2*\x)};
     \addplot[
        fill=blue,
        fill opacity=0.2
    ]
    fill between[
        of=zero and hyp,
        soft clip={domain=0.7071:2.8284},
    ];
    \addplot[name path=ellip, domain= 4/sqrt(2):4, thick, black] {sqrt(1-\x^2/16)/4};
     \addplot[
        fill=FUgreen,
        fill opacity=0.2
    ]
    fill between[
        of=zero and ellip,
        soft clip={domain=2.8284:4},
    ];
    
	\draw[black, very thin] (axis cs: 0.7071, 0.0)--(axis cs: 0.7071, 0.7071);
	\draw[black, very thin] (axis cs: 2.8284, 0.0)--(axis cs: 2.8284,  0.1768);
	
	\node[] at (axis cs: 0.35, 0.4) {\color{FUred} $\set{A}_1$};
	\node[] at (axis cs: 1.75, 0.125) {\color{blue} $\set{A}_2$};
	\draw[FUgreen!50!black] (axis cs: 3.5, 0.05)--(axis cs: 3.9, 0.3);
	\node[right=-2pt] at (axis cs: 3.9, 0.3) {\color{FUgreen!50!black} $\set{A}_3$};
	
	\draw (axis cs: 0.5, 0.8660)--(axis cs: 1, 0.9);
	\node[right=-2pt] at (axis cs: 1, 0.9) {$\ep\sqrt{1-x^2}$};
	\draw (axis cs: 2, 0.25)--(axis cs: 2.2, 0.6);
	\node[right=-2pt] at (axis cs: 2.2, 0.6) {$\displaystyle{\frac{\ep^2}{2x}}$};
	\draw (axis cs: 3.0, 0.1654)--(axis cs: 3.1, 0.8);
	\node[right=-2pt] at (axis cs: 3.1, 0.8) {$\frac{\ep}{\sigma(T)}\sqrt{1-\frac{x^2}{\ep^2\sigma(T)^2}}$};

\end{axis} 

\end{tikzpicture}

\caption{Partitioning the union of linearized pullbacks of the $\ep$-circle in the idealized situation where the right singular vectors of the fundamental matrix are not changing in time.}
\label{fig:galaxy_area}
\end{figure*}
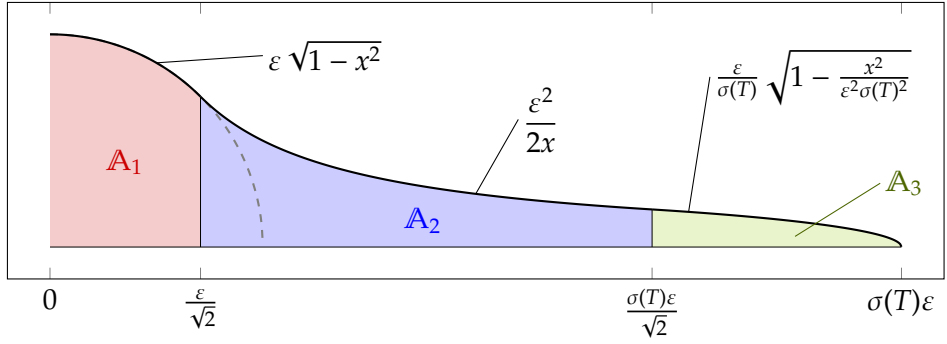

If~$\sigma(t), \sigma(t)^{-1}$ are the singular values of~$W(0,t)$, with $\sigma(t)\ge 1$, then the ellipse~$W(0,t)^{-1}\set{B}_{\ep}(0)$ is given in the $xy$-plane by the equation
\[
\frac{x^2}{\sigma(t)^2} + \sigma(t)^2y^2 \le \ep^2\,.
\]
To obtain the boundary point~$(x,y)$ of~$\set{S}(0,T)$ for fixed~$x$, we need to maximize $y$ over the possible values~$\sigma(t)\in [1,\sigma(T)]$, because the outermost ellipse defines the boundary point of the union. A straightforward computation shows that~$y = \frac{\ep^2}{2x}$, which is realized by the ellipse with semi-major axis length~$\sigma = \sqrt{2}x$. This holds for~$|x|\in[\frac{\ep}{\sqrt{2}},\frac{\sigma(T)}{\sqrt{2}}]$. For smaller~$x$ (in magnitude), the boundary coincides with that of the circle of radius~$\ep$, for larger~$x$ it is the ellipse with semi-major axis length~$\sigma(T)$ constituting the boundary. This suggests to compute the volume of~$\set{S}(0,T)$ as sum over the three distinct intervals of $x$-values,
\[
\text{vol}(\set{S}(0,T)) = 4(V_1 + V_2 + V_3),\qquad V_i = \text{vol}\,(\set{A}_i),
\]
where, as depicted in Figure~\ref{fig:galaxy_area}, 
\begin{itemize}
\item $\set{A}_1 = \left \{(x,y)\,\big\vert\, 0\le x\le \frac{\ep}{\sqrt{2}},\,  0 \le y \le \ep\sqrt{1-x^2} \right \}$;
\item $\set{A}_2  = \left \{(x,y)\,\big\vert\, \frac{\ep}{\sqrt{2}} \le x \le \frac{\sigma(T)\ep}{\sqrt{2}},\, 0 \le y \le \frac{\ep^2}{2x} \right \}$; and
\item $\set{A}_3  = \left \{(x,y)\,\big\vert\, \frac{\sigma(T)\ep}{\sqrt{2}} \le x \le \sigma(T)\ep,\, \right.$\\
$\phantom{\set{A}_3  = \left \{(x,y)\,\right\} } \left. 0 \le y \le \frac{\ep}{\sigma(T)} \sqrt{1-\frac{x^2}{\ep^2\sigma(T)^2}} \right\}$.
\end{itemize}
Direct computation gives\footnote{See {\scriptsize \href{http://www.wolframalpha.com/input/?i=integrate+sqrt(a\%5E2-x\%5E2)+from+x\%3D0+to+a\%2Fsqrt(2)}{http://www.wolframalpha.com/input/?i=integrate+sqrt(a\%5E2-x\%5E2)+from+x\%3D0+to+a\%2Fsqrt(2)}} for $V_1$. Note that $V_3 = \sigma(T)(\frac{\pi}{4}-V_1)$, since it is congruent with the patch given by the intersection of the quarter disc and $A_2$ (i.e., in Figure~\ref{fig:galaxy_area} the part of the blue shaded region under the dashed curve) stretched in the $x$ and contracted in the $y$ direction by~$\sigma(T)$, respectively. Since the stretching and contracting are by the same factor, the area of the patch under this transformation is unchanged.}
\[
V_1 = \frac{\pi+2}{8}\ep^2, \quad V_2 = \log(\sigma(T)) \frac{\ep^2}{2}, \quad V_3 = \frac{\pi-2}{8}\ep^2\,.
\]
Thus, we obtain 
\begin{proposition} \label{prop:vol_galaxy}
Under Assumption~\ref{ass:alignment} one has
\[
\text{vol}(\set{S}(0,T)) = \big(\pi + 2 \log \left( \sigma \left( T\right)\right)\big)\ep^2\,.
\]
\end{proposition}
That means, in a linear regime, where the node degree approximates $\text{vol}(\set{S}(0,T))$, one can expect an affine-linear relationship between node degree and the FTLE $\Lambda(\cdot, 0, T)=\frac{1}{T}\log \left( \sigma \left( T\right)\right)$, at least when the singular vectors do not change much in time. $\text{vol}(\set{S}(0,T)) \gg \big(\pi + 2 \log \left( \sigma \left( T\right)\right)\big)\ep^2$ is to be expected, when $\sigma \left( T\right) \gg 1$ and the singular vectors are strongly time-dependent. 

\paragraph{Singular vectors with changing direction.}
\label{sssec:changing_direction}

Assumption~\ref{ass:alignment} is not realistic in practice. However, in the following we will show reasons that it is also unrealistic that it is strongly violated in a \emph{quantitative sense}, as this would require unlikely strong vector fields.

Let $W(0,t) = U(t)\Sigma(t) V(t)^T$ be the singular value decomposition of the fundamental matrix. As above, none of its singular values are zero, and thus~$W(0,t)^{-1} = V(t) \Sigma(t)^{-1} U(t)^T$, yielding that the semiaxes of $W(0,t)^{-1} \set{B}_{\ep}(0)$ are aligned with the columns of $V(t)$, i.e., the right singular vectors of~$W(0,t)$.

Let us consider now the conditions on the dynamics that are necessary for the right singular vectors to change. As we are interested in the action of~$W(0,t)^{-1}$ on $\set{B}_{\ep}(0)$, and~$U(t)^T \set{B}_{\ep}(0) = \set{B}_{\ep}(0)$ (orthogonal transformations keep the unit ball unchanged), we set~$U(t) = I$ for simplicity. Also, by continuity of~$t\mapsto V(t)$ we have~$\det(V(t)) = 1$, and in two dimensions it means that~$V(t)$ is a rotation.

The main insight can be seen in the following prototypical example, where we rotate the columns of~$V(t)$ with angular frequency $\omega(t)$, such that
\[
\dot{V}(t) = \omega(t)\begin{pmatrix}
0 & -1\\ 1 & 0
\end{pmatrix}V(t)\,.
\]
Further, let~$\sigma_1 = \sigma_1(t)$ and $\sigma_2 = \sigma_2(t) = \frac{1}{\sigma_1(t)}$ be the diagonal elements of $\Sigma(t)$. Thus, $W(0,t) = \Sigma(t) V(t)^T$. Under our assumptions, we obtain
\begin{align*}
\dot{W} (0,t) &= \begin{pmatrix}
\dot{\sigma}_1 & 0 \\ 0 & \dot{\sigma}_2
\end{pmatrix} V(t)^T + \begin{pmatrix}
\sigma_1 & 0 \\ 0 & \sigma_2
\end{pmatrix}
\begin{pmatrix}
0 & \omega\\ -\omega & 0
\end{pmatrix}V(t)^T \\
&= \underbrace{ \begin{pmatrix}
\dot{\sigma}_1/\sigma_1 & \omega\sigma_1/\sigma_2 \\
-\omega\sigma_2 / \sigma_1 & \dot{\sigma}_2 / \sigma_2
\end{pmatrix} }_{ = D_x v, \text{ cf.~\eqref{eq:vareq}} }
W(0,t)\,,
\end{align*}
so the spatial derivative of the vector field $v$ has a component~$\omega\sigma_1/\sigma_2 = \omega\sigma_1^2$.
That is, for the right singular vector of the fundamental matrix to change its direction with unit speed the vector field $v$ needs to have a large spatial derivative, provided~$\sigma_1(t)$, the accumulated stretching from time~$t_0$ to time~$t$ is large. For hyperbolic trajectories $x(t)$ the singular value $\sigma_1(t)$ of the fundamental matrix grows exponentially in time. Thus, either the spatial derivative of the velocity field grows exponentially as well, or the change in the direction of the singular vector is exponentially slow. 

We conclude that the larger the already present local stretching in the system is, the more unlikely it is that smooth vector fields change the direction of the corresponding right singular vector significantly. With this, even if Assumption~\ref{ass:alignment} is violated, it is likely that the direction of the dominant singular vector (for trajectories showing considerable hyperbolic behavior) shows a step-function like behavior in time, as this direction is likely to change only in time intervals where~$\sigma_1(t) \ngg \sigma_2(t)$. By superposing the corresponding linearized pullback of the unit circle, one obtains a superposition of single ``galaxies'', each as in Figure~\ref{fig:ellipseunion} with different sizes and major axis directions. This is what we often observe in the examples below.
In summary, the degree correlates largely with the FTLE but, in addition, it takes rotation into account as well as nonlinear effects. In particular, the FTLE measures only stretching in the dominant direction, whereas the degree also captures expansion in the other directions. Moreover, if the singular vectors change their directions significantly, this may lead to a moderate degree even if the FTLE is small.

These theoretical considerations are underlined by the numerical example in section~\ref{ssec:DG} below.

\begin{remark}
We note that the only results that require the dynamics to be area-preserving are Proposition~\ref{prop:vol_galaxy}, the considerations is Section~\ref{sec:degree}b, and those in Figure~\ref{fig:overlap}.
These results rely on the fact that the fundamental matrix of the system has two singular values whose product is one.
The qualitative classification of dynamical regimes in Table~\ref{tab:regimes} is independent of area-preservingness.
\end{remark}

\subsection{Clustering coefficient.}
\label{ssec:clco}

Recall that the clustering coefficient \eqref{eq:cc} of a node $x_i$ is defined as
\begin{equation}
c_i = \frac{\sum_{j,k}^N A_{ij} A_{jk}A_{ki}}{d_i(d_i - 1)}
\end{equation}
where $A_{ij}$ is as in \eqref{eq:A}. It counts the number of triangles with vertex $x_i$ divided by the total number of possible triangles. In the limit of very many data points (and thus large degree), we can approximate the clustering coefficient by simplifying the denominator above, yielding
\begin{equation}
\tilde c_i = \frac{\sum_{j,k}^N A_{ij} A_{jk}A_{ki}}{d_i^2}.
\end{equation}
From now on this equation will be used, and the tilde will be dropped. It is useful to write this as
\[
c_i = \frac{1}{d_i} \sum_j A_{ij} f_j^{(i)}, \qquad f_j^{(i)} = \frac{1}{d_i}\sum_k A_{jk}A_{ki}.
\]

\paragraph{Geometrical interpretation of the clustering coefficient.}
We already know that $N^{-1}d_i$ converges to $\text{vol}\left(\set{G}_{\ep} (x_i)\right)$, i.e., the volume of the galaxy neighborhood, in the infinite data limit $N\rightarrow \infty$. Further, we may write $A_{ij} = \mathbf{1}_{\set{G}_{\ep}(x_i)}(x_j)$ and then $f_j^{(i)}$ converges in the Monte Carlo sense to the function $f^{(i)}$ evaluated at the point $x_j$:
\begin{align*}
f^{(i)}(x_j) &\; = \frac{1}{\text{vol}(\set{G}_{\ep}(x_i))}\int \mathbf{1}_{\set{G}_{\ep} (x_j)}(y) \mathbf{1}_{\set{G}_{\ep}(x_i)}(y) dy \\
&\; = \frac{1}{\text{vol}(\set{G}_{\ep}(x_i))}\int \mathbf{1}_{\set{G}_{\ep}(x_j) \cap \set{G}_{\ep}(x_i)}(y) dy \\
&\; = \frac{{\text{vol}(\set{G}_{\ep}(x_i)\cap \set{G}_{\ep}(x_j))}}{{\text{vol}(\set{G}_{\ep}(x_i))}}.
\end{align*}
This is the proportion of $\set{G}_{\ep}(x_i)$ that overlaps $\set{G}_{\ep}(x_j)$.

Finally, the clustering coefficient converges to
\begin{equation}
\label{eq:cc_estim}
\begin{aligned}
c(x_i) &\; =  \frac{1}{\text{vol}(\set{G}_{\ep}(x_i))} \int \mathbf{1}_{\set{G}_{\ep}(x_i)}(y) f^{(i)}(y) dy \\
&\; =  \frac{1}{\text{vol}(\set{G}_{\ep}(x_i))} \int_{\set{G}_{\ep}(x_i)}  \frac{{\text{vol}(\set{G}_{\ep}(x_i)\cap \set{G}_{\ep}(y))}}{{\text{vol}(\set{G}_{\ep}(x_i))}} dy.
\end{aligned}
\end{equation}
In other words, $c(x)$ is the expected relative overlap of the neighborhood $\set{G}_{\ep}(x)$ and a second neighborhood $\set{G}_{\ep}(y)$ where $y$ is 
drawn from $\set{G}_{\ep}(x)$. 

In the linear regime, that is, for small $\ep$ and appropriate time spans, these neighborhoods and their volumes can be approximated using the variational equation and the estimates in the preceding subsection.
When both $\set{G}_{\ep}(x)$ and $\set{G}_{\ep}(y)$ are balls the expected overlap can be explicitly computed, but this is already no longer possible when ellipsoids have to be taken into account. 

\begin{figure}[h]
\centering
\includegraphics[width = 0.49\textwidth]{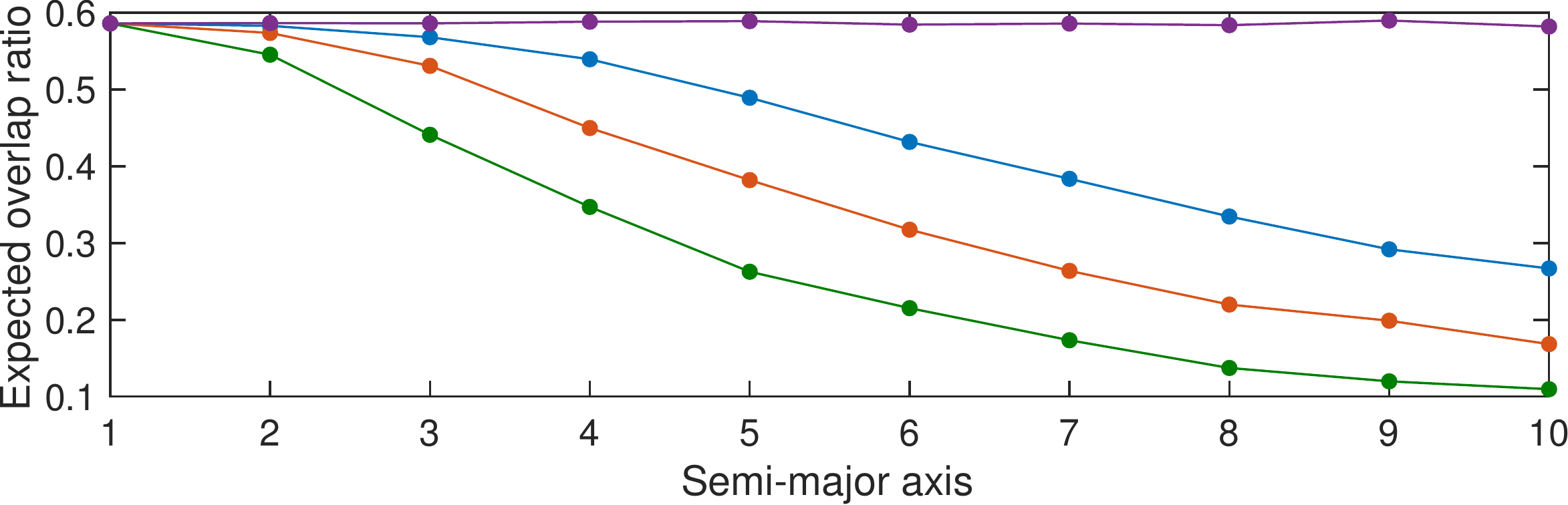}
\caption{Expected relative overlaps of two ellipses depending on the semi-major axis and the rotation angle ($0$ (purple),  up to $\pm \frac{\pi}{32}$ (blue), $\pm \frac{\pi}{16}$ (red), $\pm \frac{\pi}{8}$ (green)).}
\label{fig:overlap}
\end{figure}

In Figure \ref{fig:overlap}, we have numerically (via a Monte Carlo approach) estimated the expected relative overlaps of equally sized ellipses in the 2D case depending on the length of the semi-major axis $\sigma\geq 1$, with the semi-minor axis being $\sigma^{-1}$. The purple curve corresponds to the (unrealistic) case, when the overlapping ellipses are exactly axis-parallel. Here, the numerical estimation is independent of $\sigma$ and nicely matches the theoretical value for two overlapping circles of $1-\frac{3\sqrt{3}}{4\pi}\approx 0.59$.\footnote{The integral $$ \frac{1}{\pi}\int_{0}^1 4x\arctan{\frac{x}{2}}-x^2\sqrt{4-x^2} \, dx= 1-\frac{3\sqrt{3}}{4\pi}$$ exactly describes the expected relative overlap of two circles $\set{B}_{\ep}(x)$ and $\set{B}_{\ep}(y)$, where $y$ is uniformly drawn from $\set{B}_{\ep}(x)$.} 
For the other curves, the overlapping ellipses are allowed to be slightly rotated (rotation angles up to $\pm \frac{\pi}{32}$ (blue), $\pm \frac{\pi}{16}$ (red), $\pm \frac{\pi}{8}$ (green)) and in these cases the expected relative overlap decreases as $\sigma$ is increased. In particular, an increase in the maximum rotation angle also leads to a decrease in the expected overlap ratio when the length of the semi-major axis are kept fixed.

This confirms the frequent observation that the local clustering coefficient takes large values where the dynamics is elliptic\cite{Donner_et_al_2010a,RoSeHe17}. In this case, the corresponding galaxy neighborhood is ball-like and thus the FTLE and also the degree are small. Moreover, the local clustering coefficient is small when there is strong stretching and the FTLE is large. So we expect that the FTLE and the local clustering coefficient are strongly negatively correlated. However, due to the finite $\ep$ in the network construction the local clustering coefficient measures also nonlinear effects as we will discuss in the following.

\paragraph{Time-dependent behavior.} 
Note that by construction, the set of neighbors in the network increases in a nested manner as the observation time interval increases: If $A_{ij} = 1$ for trajectories observed for the time interval $[t_0,t]$, naturally $A_{ij} = 1$ holds also for the time interval~$[t_0,t']$ with~$t' > t$.
Unlike the degree, which thus increases monotonically in time,
we expect the \emph{qualitative} behavior of the clustering coefficient to change in time considerably.

In a dynamically mixing region, the clustering coefficient starts for small times with a moderate to large value, as the network for small times is based on vicinity of initial points. Then, as the observation time grows, it decreases, as new neighbors are introduced which are not necessarily neighbors of neighbors (due to hyperbolic stretching). This holds in an intermediate time interval for which the image of small balls under the dynamics is a filamented set, but does not yet fold back to itself and ``cover'' full-dimensional subsets of the state space. Then, as time increases, we expect the clustering coefficient to grow again, as eventually any two points get close-by in a mixing region again\footnote{For dynamics that is mixing in the measure-theoretic sense, this is shown, e.g., in Footnote~5 of Koltai \& Renger\cite{KoRe18}.}, and the filaments tend to become ``space-filling''.
To be more precise, once the filamentation is so strong that in a measure-theoretic sense $\set{G}_{\ep}(x_i)\cap \set{G}_{\ep}(y) \approx \set{G}_{\ep}(x_i)$ for every~$y\in \set{G}_{\ep}(x_i)$, then by~\eqref{eq:cc_estim} the clustering coefficient gets large again.

In regular regions---where the dynamics is not distorting strongly and the mutual distances between points barely change---we expect the clustering coefficient to stay approximately constant all the time.

We summarize the expected characteristics of the network measures which we expect for different finite-time dynamical behavior in Table~\ref{tab:regimes}.
\begin{table}[h]
\begin{tabular}{l | l | l}
\textbf{degree} & \begin{tabular}[t]{@{}l@{}} \textbf{clustering}\\ \textbf{coefficient}\end{tabular} & \begin{tabular}[t]{@{}l@{}} \textbf{finite-time dynamical}\\ \textbf{regime} \end{tabular} \\ \hline
\begin{tabular}[t]{@{}l@{}} small --\\ moderate \end{tabular}  & large & elliptic or parabolic motion \\ \hline
moderate & small & \begin{tabular}[t]{@{}l@{}} filamentation (finite-time\\ hyperbolicity), no mixing \end{tabular} \\ \hline
moderate & moderate & \begin{tabular}[t]{@{}l@{}} stickiness (mixing \\ close to regular regions)   \end{tabular} \\ \hline
large & large & mixing
\end{tabular}
\caption{Network measures and finite-time dynamical regimes.}
\label{tab:regimes}
\end{table}

\section{Discussion of further network measures.}
\label{sec:furthermeas}
Let us now briefly discuss the expected behavior of the other network measures from above. To this end it is helpful to differentiate two kinds of dynamical behaviors that are very characteristic of complicated flows we are interested in.
The first we connect to ``mixing regions'', where (weak) mixing is understood in the measure-theoretic sense~\cite{Wal00}.
The second kind is connected to ``regular regions'', and refer to those regions in state space that are not mixing, and we think of them as regions performing a rigid-body motion---up to slight distortions. Naturally, this is not a complete or well-defined partition of the flow domain, as it uses notions (like mixing) that are defined for infinite time, and we are looking at finite time intervals;  so there could be a whole homotopy of characteristics connecting these two. Still, as a descriptor of the two ``extreme cases'', it will prove very useful in the following.

\paragraph{Closeness.}

As network measures get more complicated, it gets significantly more involved to connect them directly to the dynamical behavior of a system.
It is safe to claim that in a mixing region most trajectories eventually ``meet'' one another, their mutual distances will be small, giving a large closeness value. In a regular region the mutual arrangement of trajectories stays similar in time, thus their distance stays moderate to large, and hence we expect them to have a moderate closeness.

Because the set of neighbor nodes increases in a nested manner in time, the length of shortest paths decrease monotonically as the observation time increases, and thus closeness increases. In contrast to the clustering coefficient, we do not expect closeness to change its qualitative behavior in time: Based on the above considerations, closeness increases more rapidly in mixing regions than in regular regions, but the qualitative picture with respect to this network measure does not change with the length of the time interval in consideration.


Note that closeness here works with distances of trajectories with respect to shortest paths in the graph given by~$A$. A concept in its nature similar to this was put forward in Ref.~\onlinecite{KoRe18}, where a ``semidistance of mixing''\footnote{This semidistance is defined as the shortest path in a time-dependent graph, where the weight of edge $(i,j)$ at the $k$-th time instance is the squared distance of the trajectories $i$ and $j$ at time $t_k$, ${\| x_i(t_k) - x_j(t_k) \|^2}$, while in every step it is allowed to stay in the same node (self-transitions have zero weight). Thus, this semidistance is short between two trajectories that eventually come close during the time of consideration.} for a finite set of Lagrangian trajectory data is derived. Further, it is shown that it can be computed by shortest paths in time-dependent graphs that comprise similar information to the $\ep$-neighbor adjacency graphs at some fixed time instant~$t$ (without accumulating the neighbors in time, as done in~\eqref{eq:A}). More crucially, it is demonstrated that coherent sets---sets that move with the flow and do not mix with their exterior while doing so---are regions ``maximally far'' from one another with respect to this distance. This connection suggests that coherent regions have a large mean distance to the rest of the network, and that their closeness is thus smaller. Our observations for the one-dimensional example in section~\ref{ssec:ex1d} confirm this. As closeness is much more expensive to compute than the other network measures considered here (as shortest paths between \emph{all} pairs of node need to be computed, giving a best-case complexity of~$\mathcal{O}(n^3)$ for the Floyd--Warshall algorithm, and $\mathcal{O}(n^2\log(n))$ for the Dijkstra algorithm with Fibonacci heap), we refrain from further numerical analysis of closeness.

Finally, we note that closeness is influenced also by the relative position of a trajectory with respect to the other trajectories; towards the ``boundary'' of the trajectory ensemble closeness is expected to be smaller. This is nicely reflected in the experimental example of section~\ref{ssec:aviso}.

\paragraph{Betweenness centrality.}

The betweenness centrality $b_i$ of a node $x_i$ measures the proportion of shortest paths of the network that traverse~$x_i$. Thus---like closeness---it takes global information of the network into account, and computationally it has the same complexity. 
In general networks, $b$ is large for nodes that connect different almost-decoupled subgraphs. Such nodes, often called \emph{hubs}, play a crucial role in the transfer of information or mass transport. So, nodes with high betweenness centrality in our trajectory-based network correspond to structures that connect the different coherent regions of the underlying system. 

\paragraph{Average node degree.} 

While the local degree of a node $x_i$ only takes the direct neighboring nodes of $x_i$ into account, by taking their average degree $\langle d \rangle_{nn, i}$  an extended neighborhood is considered. The resulting average degree field is a smoothed version of the degree field, with similar properties and with $d$ and $\langle d \rangle_{nn}$ strongly correlated. Significant quantitative differences in the fields may only occur at nodes $x_i$ where the node degree is locally maximal or minimal. 
Such a difference is measured by the degree anomaly $\Delta d_i= d_i - \langle d \rangle_{nn, i}$, which may serve as an indicator of the local heterogeneity of the phase space structures\cite{Donner_et_al_2010a}.
\section{Numerical examples}
\label{sec:numerics}

\subsection{One-dimensional prototypical example}
\label{ssec:ex1d}

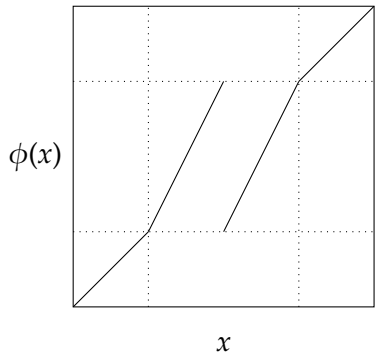
\begin{figure}[h]
\begin{tikzpicture}[scale=1]
  \tikzstyle{every node}=[font=\small]
  \draw(0,0) rectangle (4,4);
  \draw[dotted] (0,1)--(4,1);
  \draw[dotted] (0,3)--(4,3);
  \draw[dotted] (1,0)--(1,4);
  \draw[dotted] (3,0)--(3,4);
  \draw(0,0)--(1,1);
  \draw(1,1)--(2,3);
  \draw(2,1)--(3,3);
  \draw(3,3)--(4,4);
  \node at (2,-0.5) {$x$};
  \node at (-0.5,2) {$\phi(x)$};
\end{tikzpicture}

\caption{The system~\eqref{eq:1dsys}.}
\label{fig:1dsys}
\end{figure}

To gain some intuition about regular and mixing regions in a ``controlled environment'', let us now consider the discrete-time system on~$\set{M}=[0,1]$ given by
\begin{equation} \label{eq:1dsys}
{
\renewcommand{\arraystretch}{1.25}
\phi(x) = \left\{ \begin{array}{ll}
x, & x \in [0,\tfrac14) \cup (\tfrac34,1] \\
\left(2(x-\tfrac14) \mod \tfrac12\right) + \tfrac14, & x \in [\tfrac14,\tfrac34]\,,
\end{array}
\right.
}
\end{equation}
see Figure~\ref{fig:1dsys}. This map has three invariant sets. The left and right ones are static, such that the mapping restricted to them is the identity, and are meant to model regions of the state space in complicated flows, that are ``regular'' in the sense that the mutual distance of points is not changed (or just barely) by the dynamics. We will consider these as one kind of prototype for coherent vortices. The third region physically separates the other two, and the dynamics on it is mixing (as it is the circle doubling map).

\begin{figure}[h]
\centering

\includegraphics[width = 0.23\textwidth]{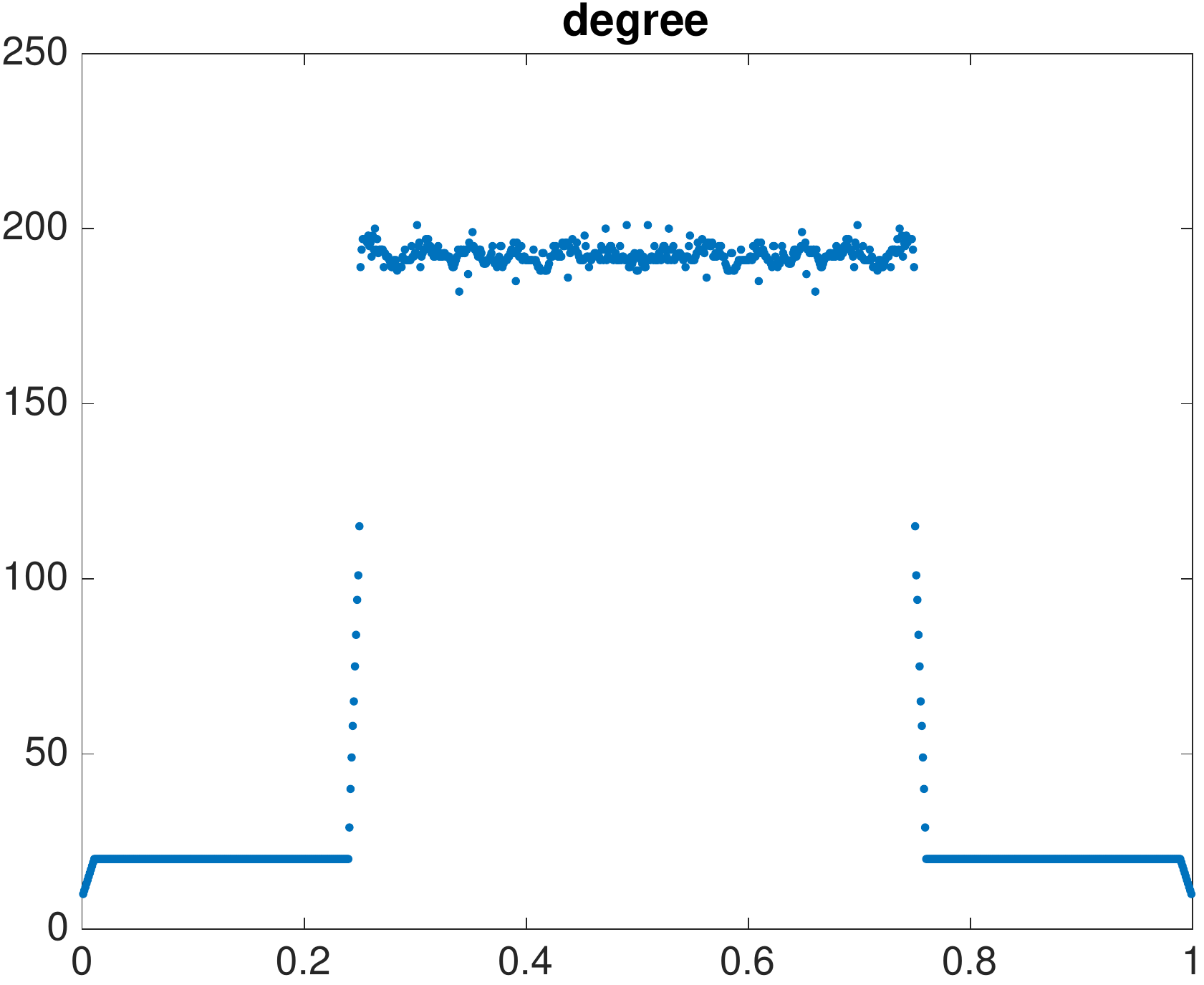}
\hfill
\includegraphics[width = 0.23\textwidth]{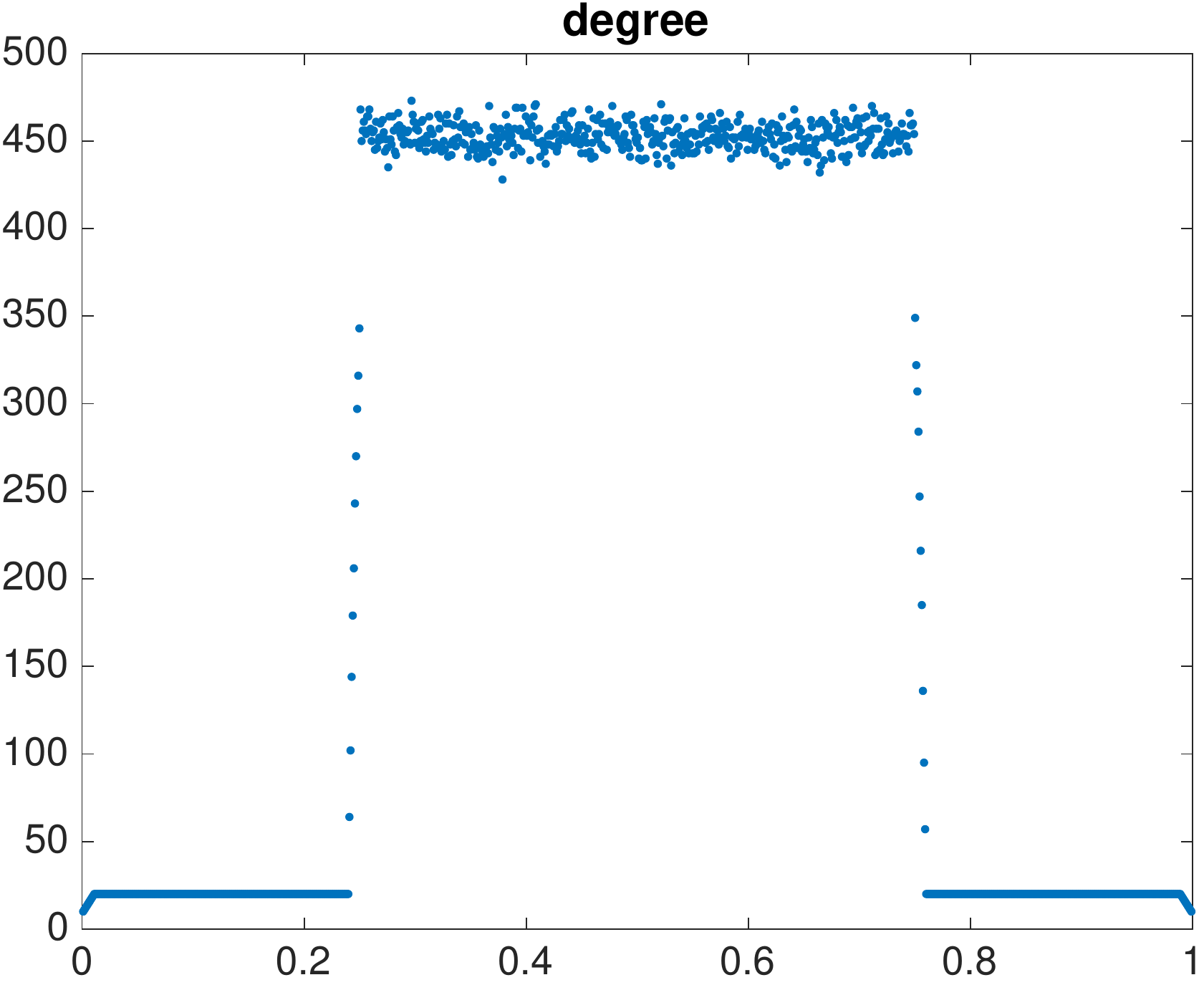}

\includegraphics[width = 0.23\textwidth]{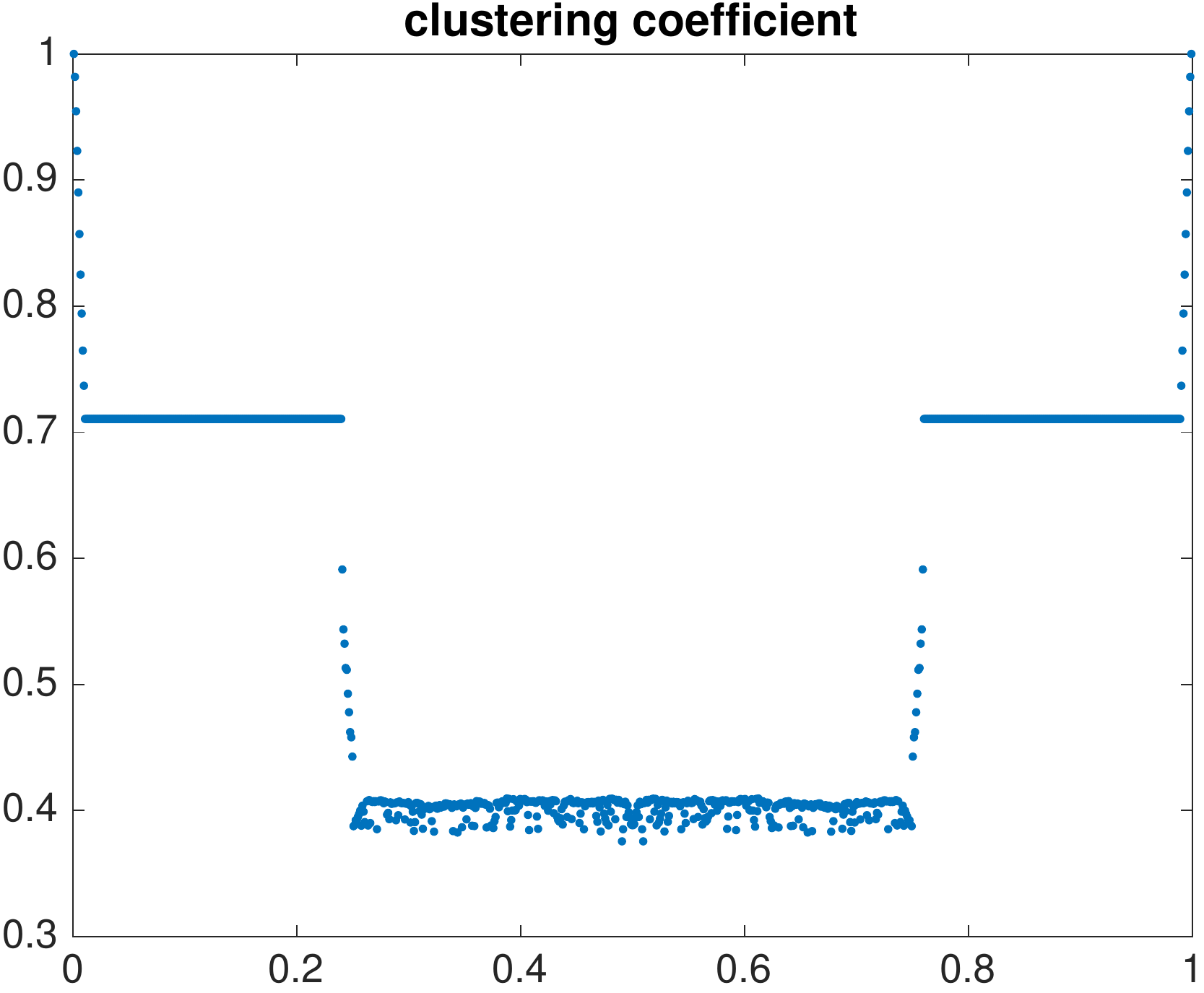}
\hfill
\includegraphics[width = 0.23\textwidth]{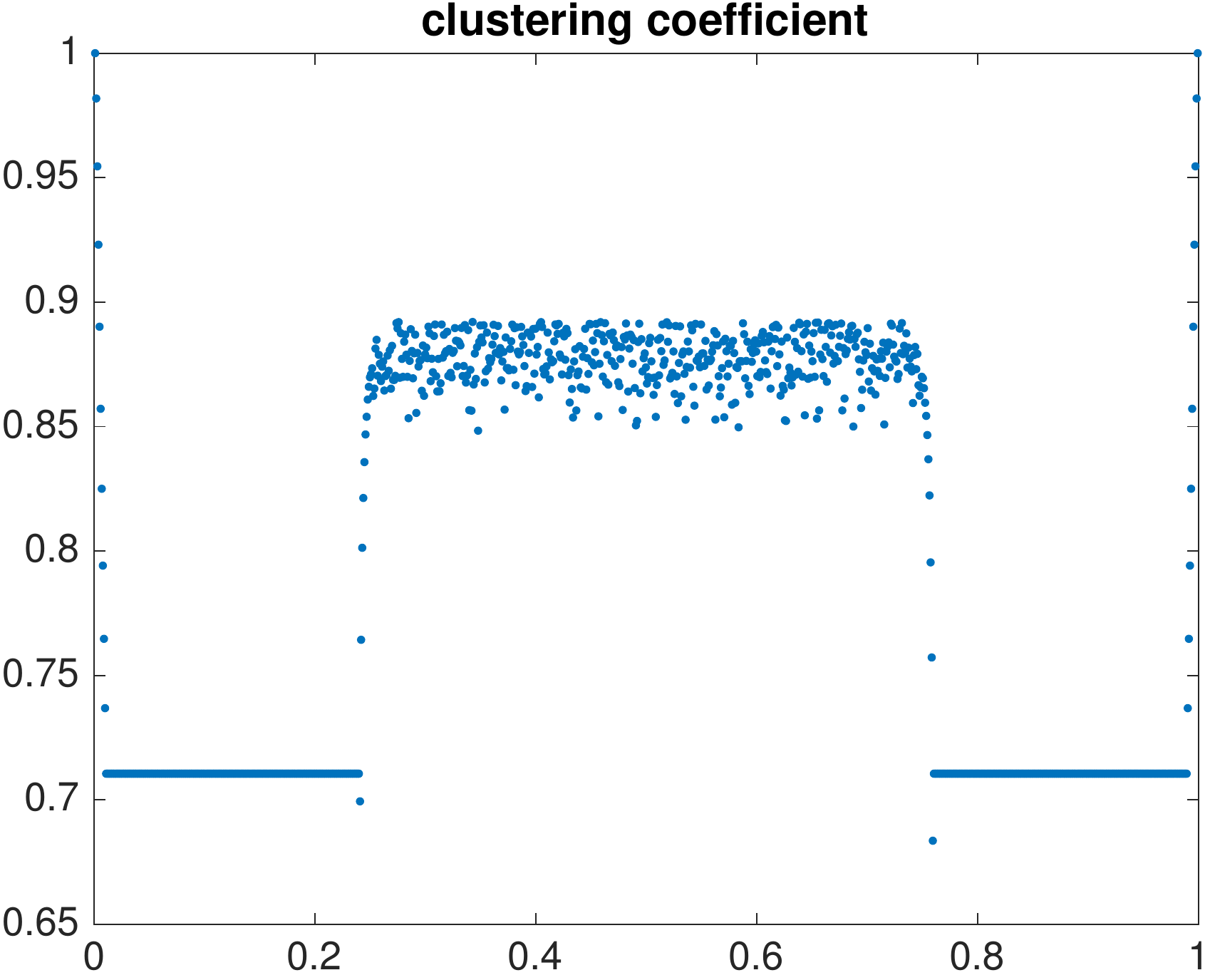}

\includegraphics[width = 0.23\textwidth]{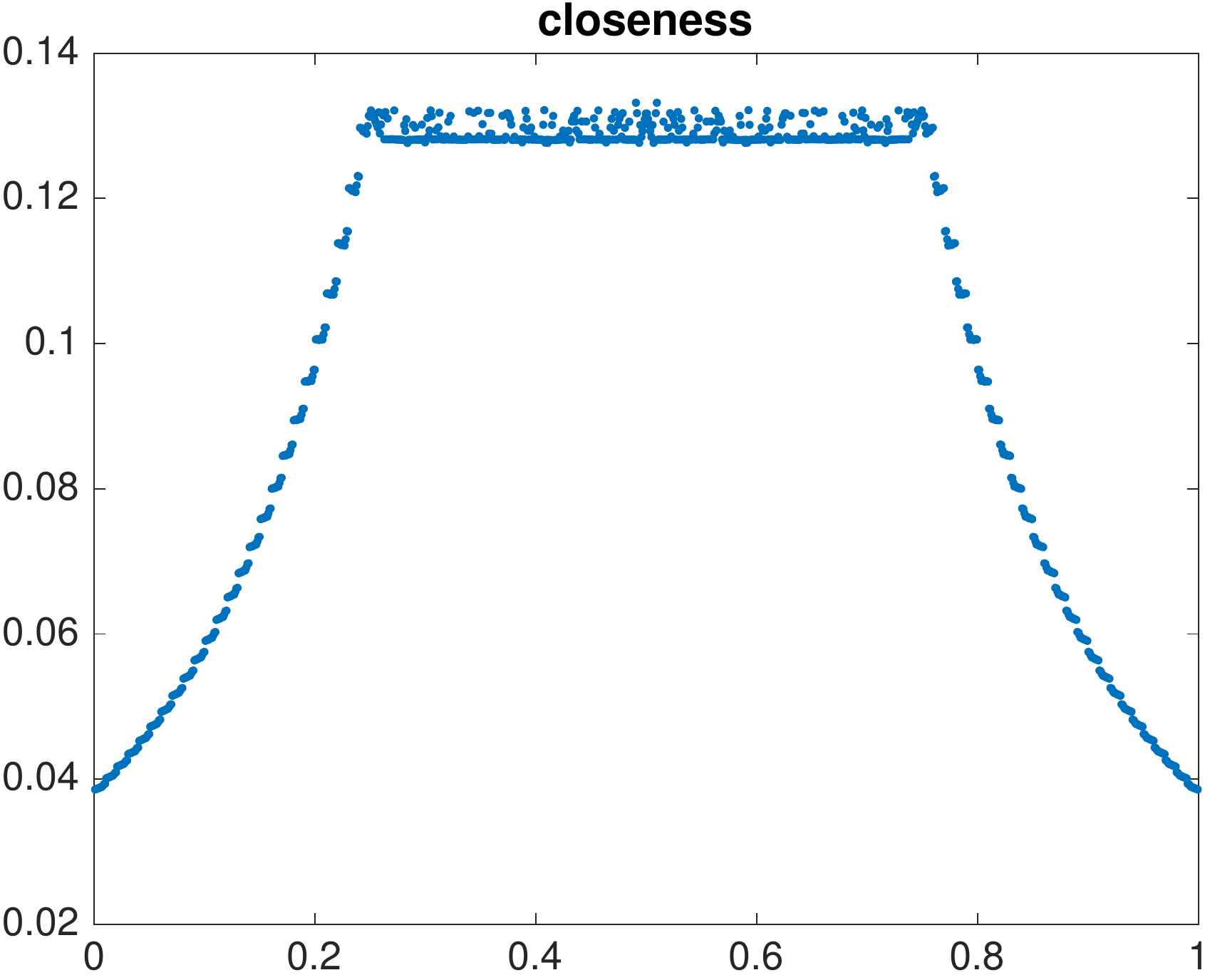}
\hfill
\includegraphics[width = 0.23\textwidth]{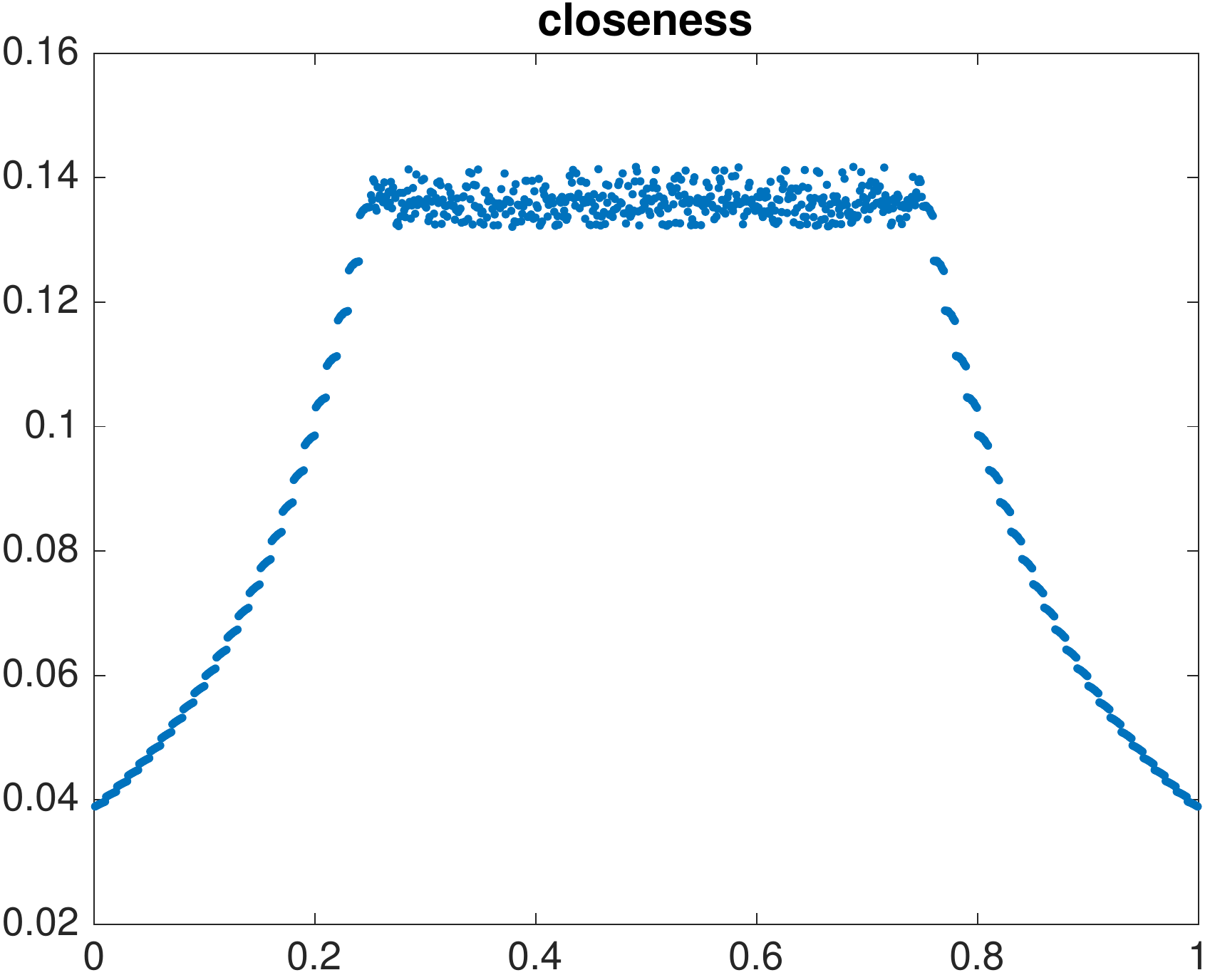}


\caption{Network measures and resulting classification for 20 (left) and 100 (right) steps of the prototypical one-dimensional system~\eqref{eq:1dsys}. Note that only the behavior of the clustering coefficient changes \emph{qualitatively}, as we also expected.
}
\label{fig:1dsys_measures}
\end{figure}

\begin{figure}[h]
\centering

\includegraphics[width = 0.23\textwidth]{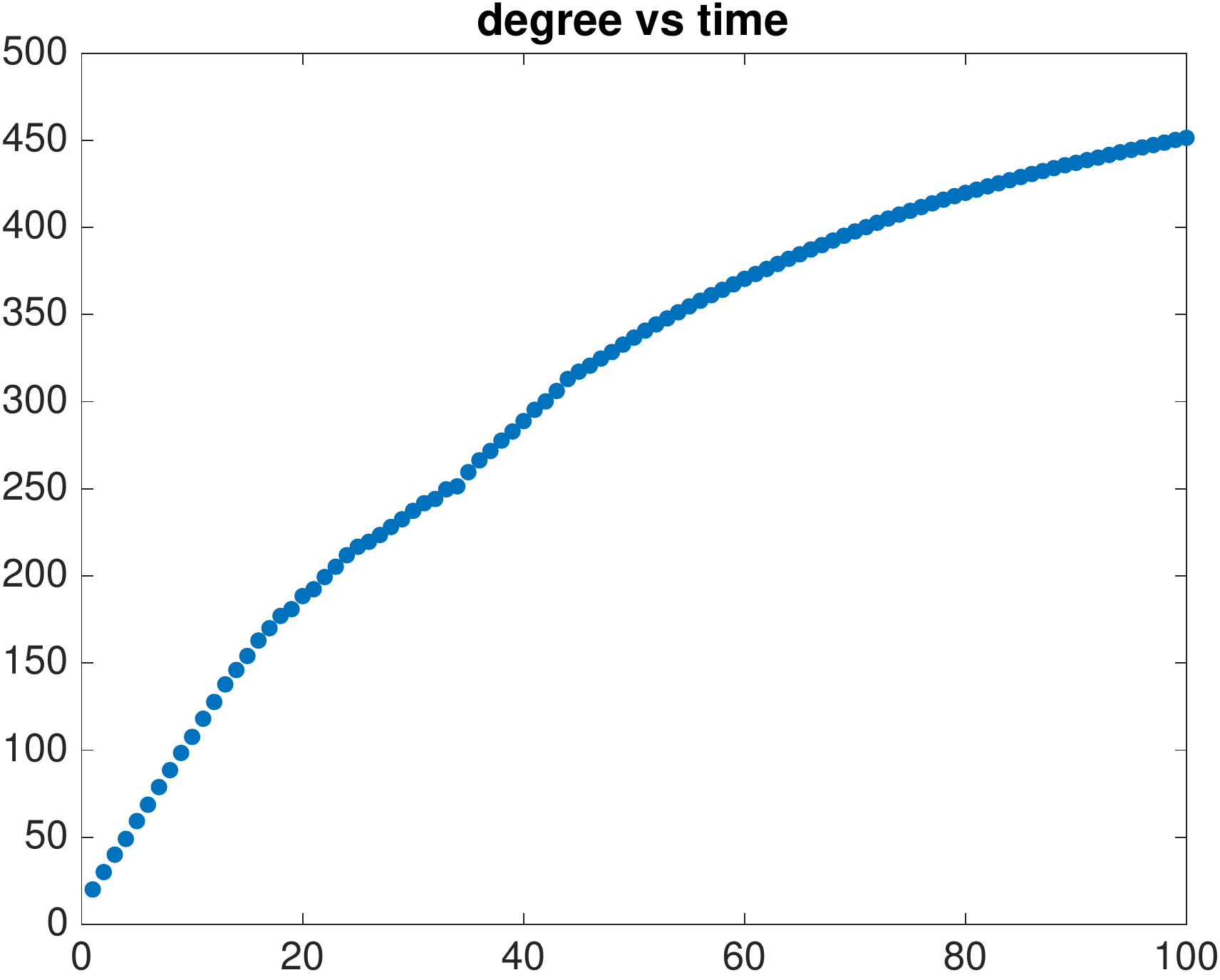}
\hfill
\includegraphics[width = 0.23\textwidth]{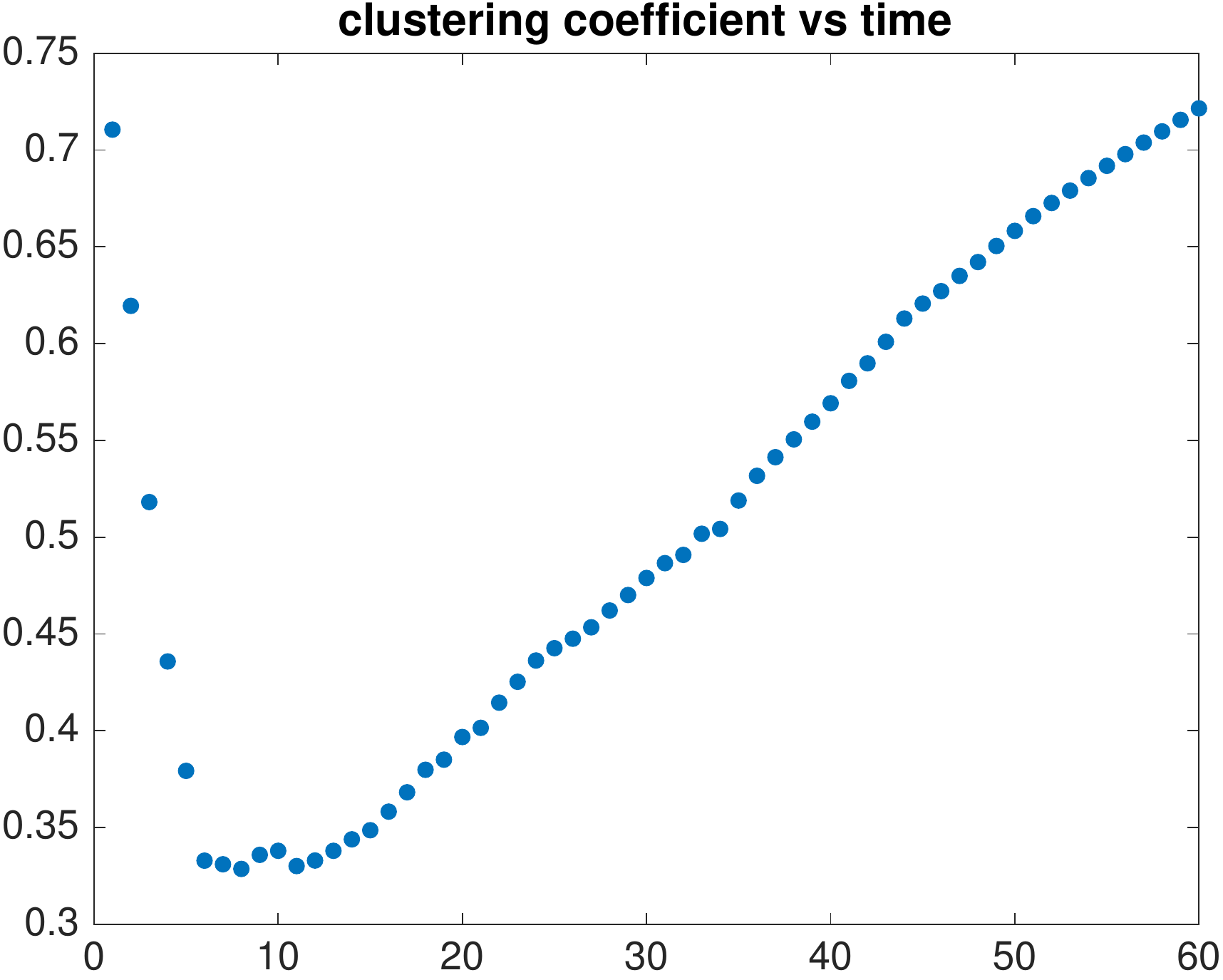}

\caption{The evolution of degree and clustering coefficient (averaged for the mixing region $\set{M}_0 = [1/4,3/4]$) as the time window $[0,t]$ grows. As the mapping on $\set{M}_0$ is the circle doubling map, in every iteration half of the neighborhood of every point is exchanged, thus the degree grows initially by $\sim 9$ per iteration. As the time window grows, this value increases, as by mixing every point eventually gets arbitrarily close to any other in~$\set{M}_0$, and eventually there are no new neighbors to find. The clustering coefficient also shows the behavior we expect; first dropping, as the neighborhood grows rapidly by mixing, then growing, as arbitrary two neighbors are likely to meet at some point, also due to mixing.}
\label{fig:1dsys_timeevolution}
\end{figure}

We carry out our computations for a network with 1000 initially equispaced trajectories, and $\ep = 0.01$. By this choice, all points (neglecting the boundaries) have initially 18 neighbors, resulting 
in $\frac{18\cdot 17}{2}$ triples. The initial number of triangles of a node $x_i$ is $108$ and thus the proportion of triples that are triangles is $\frac{2\cdot 108}{18\cdot 17}\approx 0.706$. These quantities, $18$ and $0.706$, coincide with the degree and the local clustering coefficient, respectively, for all times in the static regions, as shown in Figure~\ref{fig:1dsys_measures}.  Moreover, this study confirms what we expected from the analysis in section~\ref{sec:degree}: the degree and closeness are monotonic in the observation time span, but do not change qualitatively, while the clustering coefficient shows a very clear qualitative change between the regular and mixing region as the time span grows. This is highlighted in Figure~\ref{fig:1dsys_timeevolution}, where the time-evolution of the degree and clustering coefficient is shown, averaged over the mixing region.

\subsection{Periodically driven double gyre flow}
\label{ssec:DG}
As a benchmark problem for analyzing flow structures we consider the double gyre flow\cite{shadden_lekien_marsden_05}, a time-dependent system of differential equations 
\begin{eqnarray}
\dot{y} &=& -\pi A \sin(\pi f(y, t))\cos(\pi z) \label{eq: dgyre_eqn}\\
\dot{z} &=& \pi A \cos(\pi f(y, t))\sin(\pi z)\frac{df}{dy}(y, t),\notag
\end{eqnarray}
where $x = (y,z)\in\R^2$ is the state,
$f(y,t)=\delta \sin(\omega t)y^2+(1-2\delta\sin(\omega t))y$.
We choose parameter values $A=0.25$,
$\delta=0.25$, $\omega=2\pi$ and fix $t_0=0$. We obtain a flow of period
$\tau=1$ on the domain $\set{M}=[0,2]\times[0,1]$. 

\begin{figure*}[htbp]
\centering
\includegraphics[width = 0.33\textwidth]{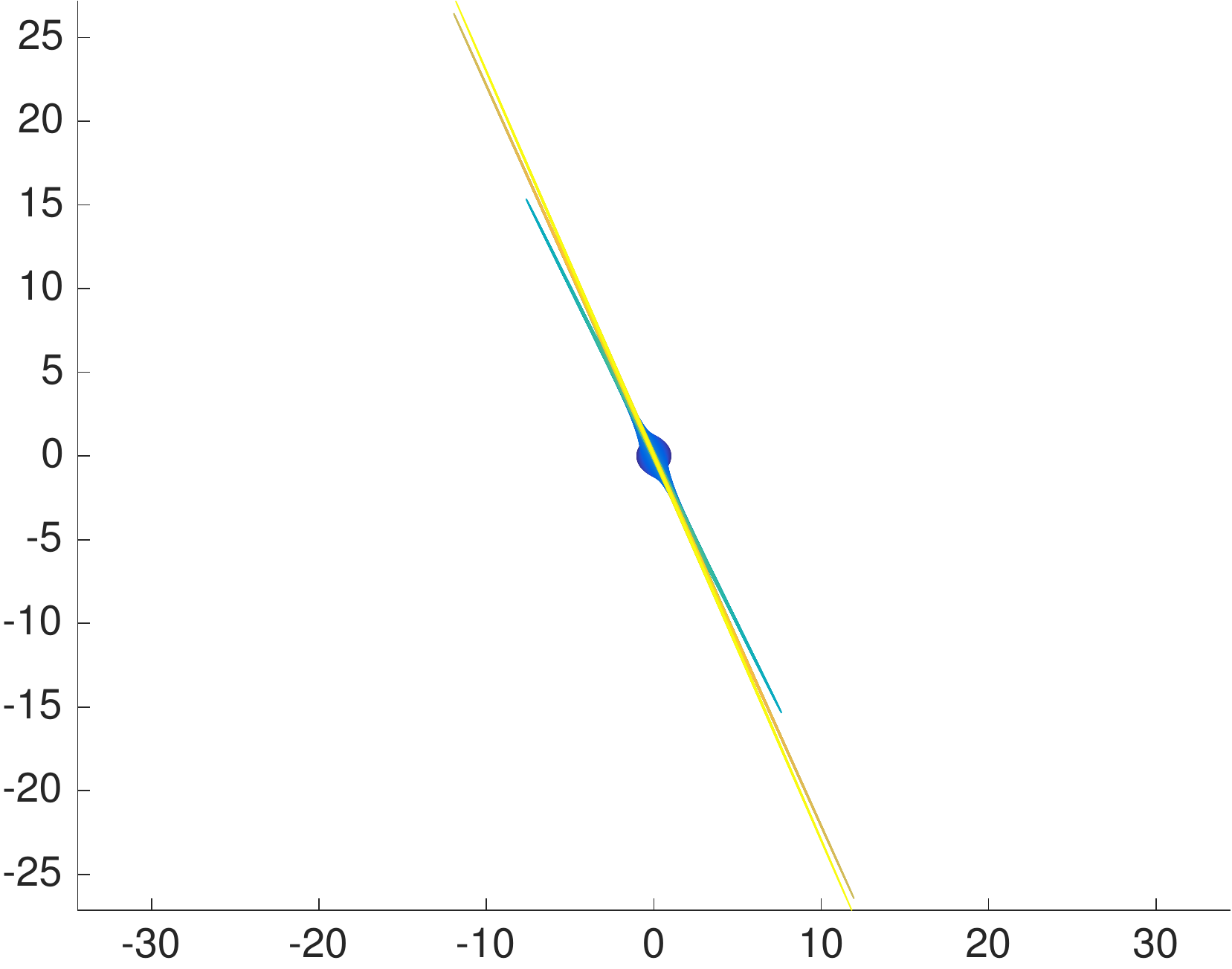}
\hfill
\includegraphics[width = 0.33\textwidth]{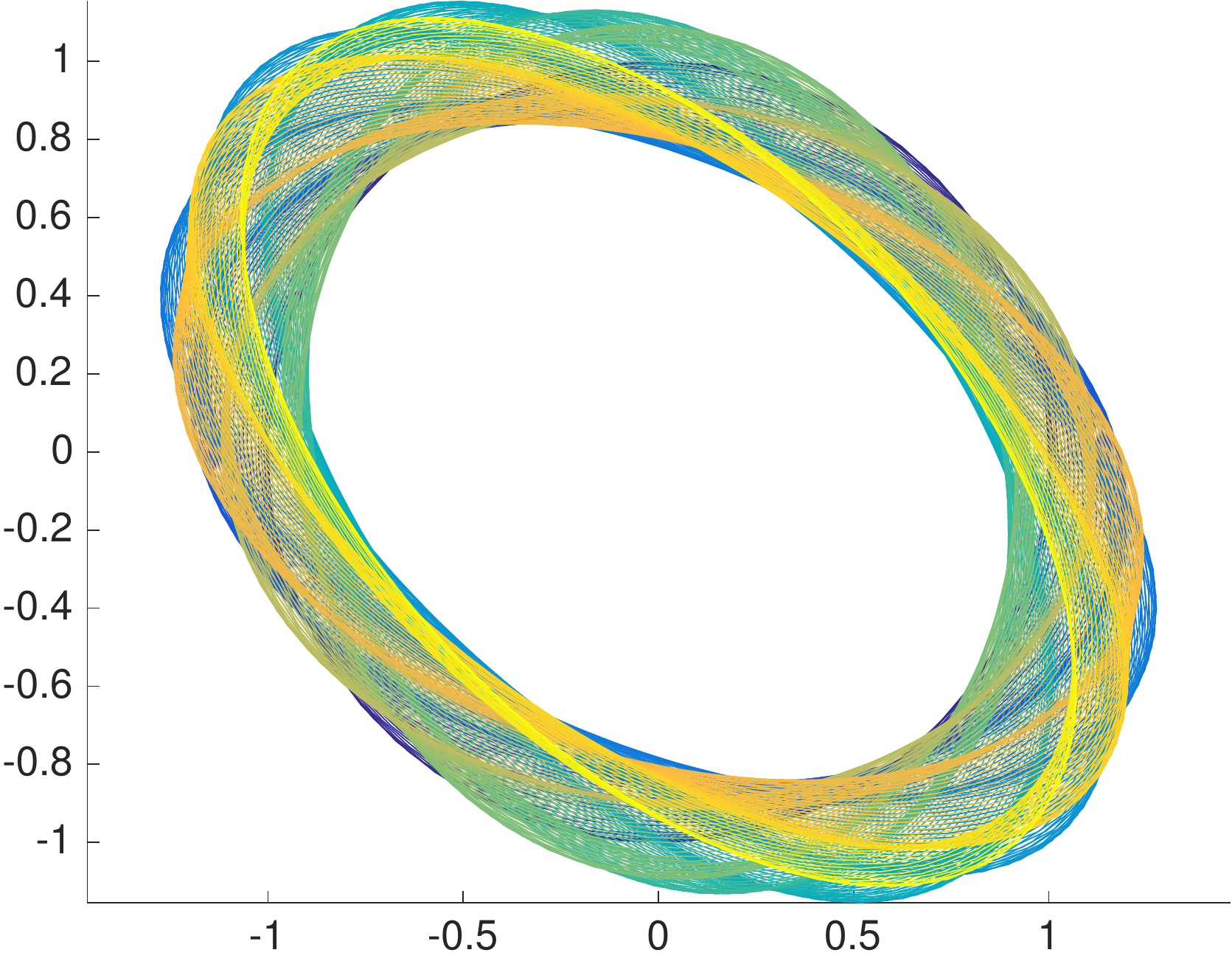}
\hfill
\includegraphics[width = 0.33\textwidth]{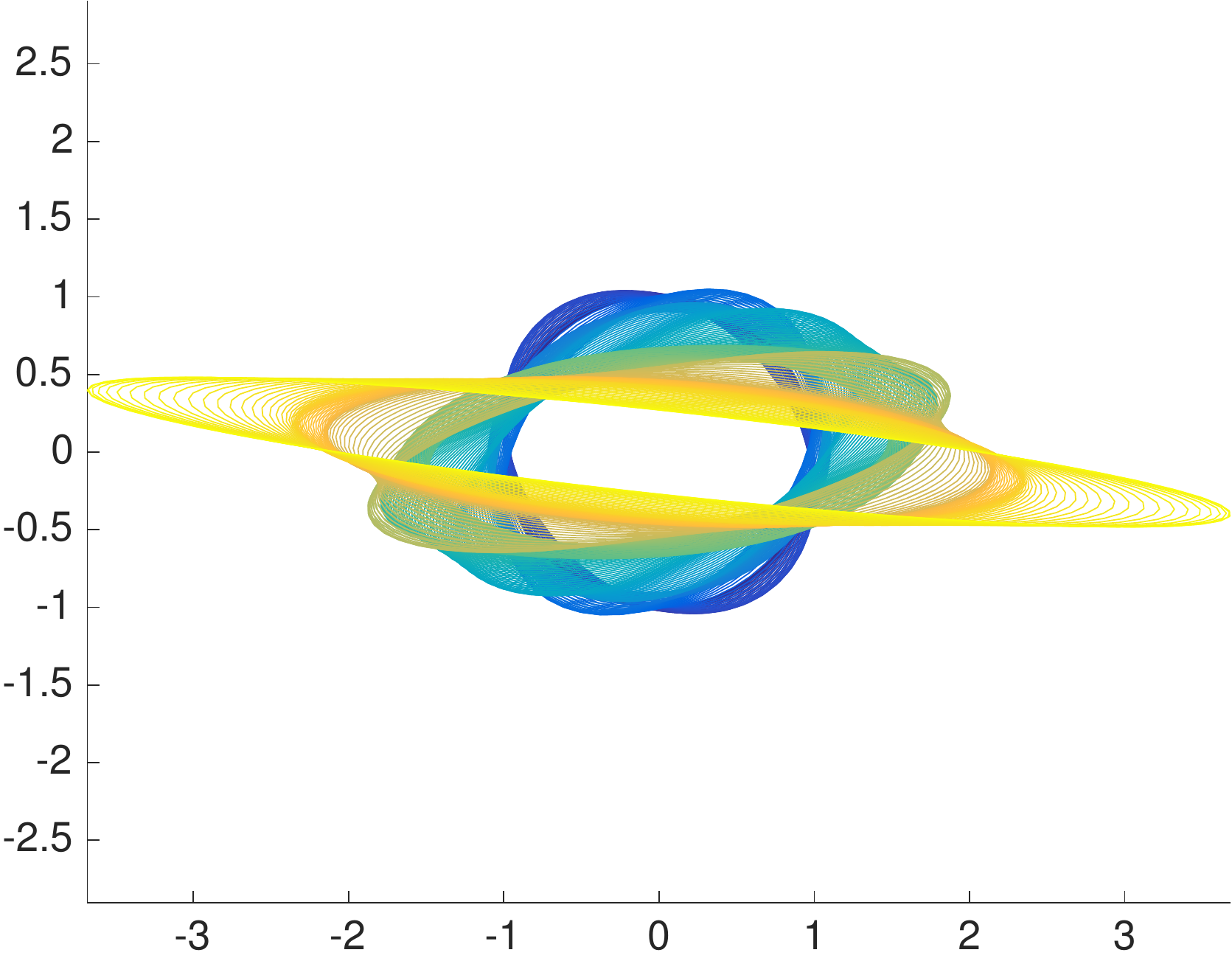}

\includegraphics[width = 0.33\textwidth]{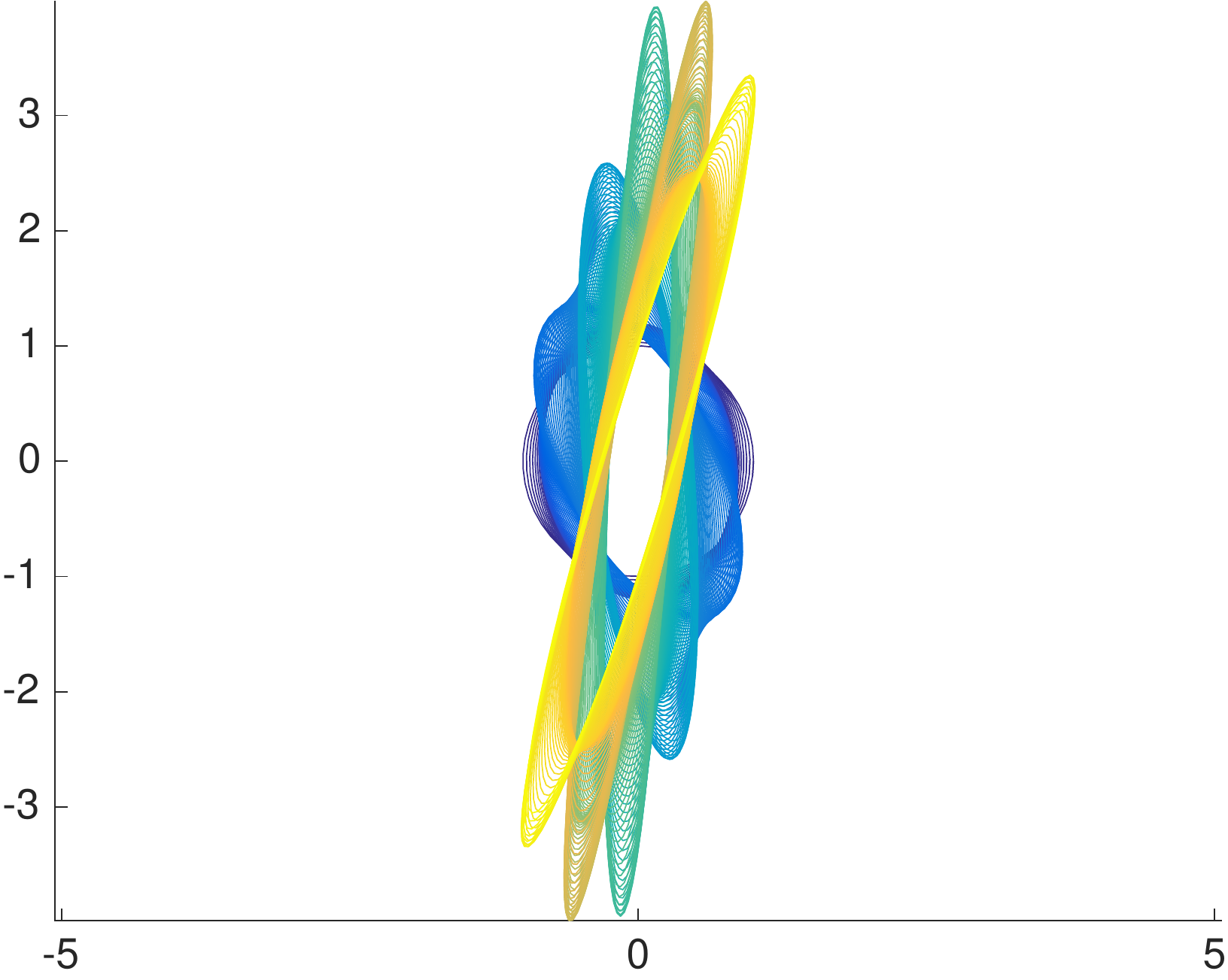}
\hfill
\includegraphics[width = 0.33\textwidth]{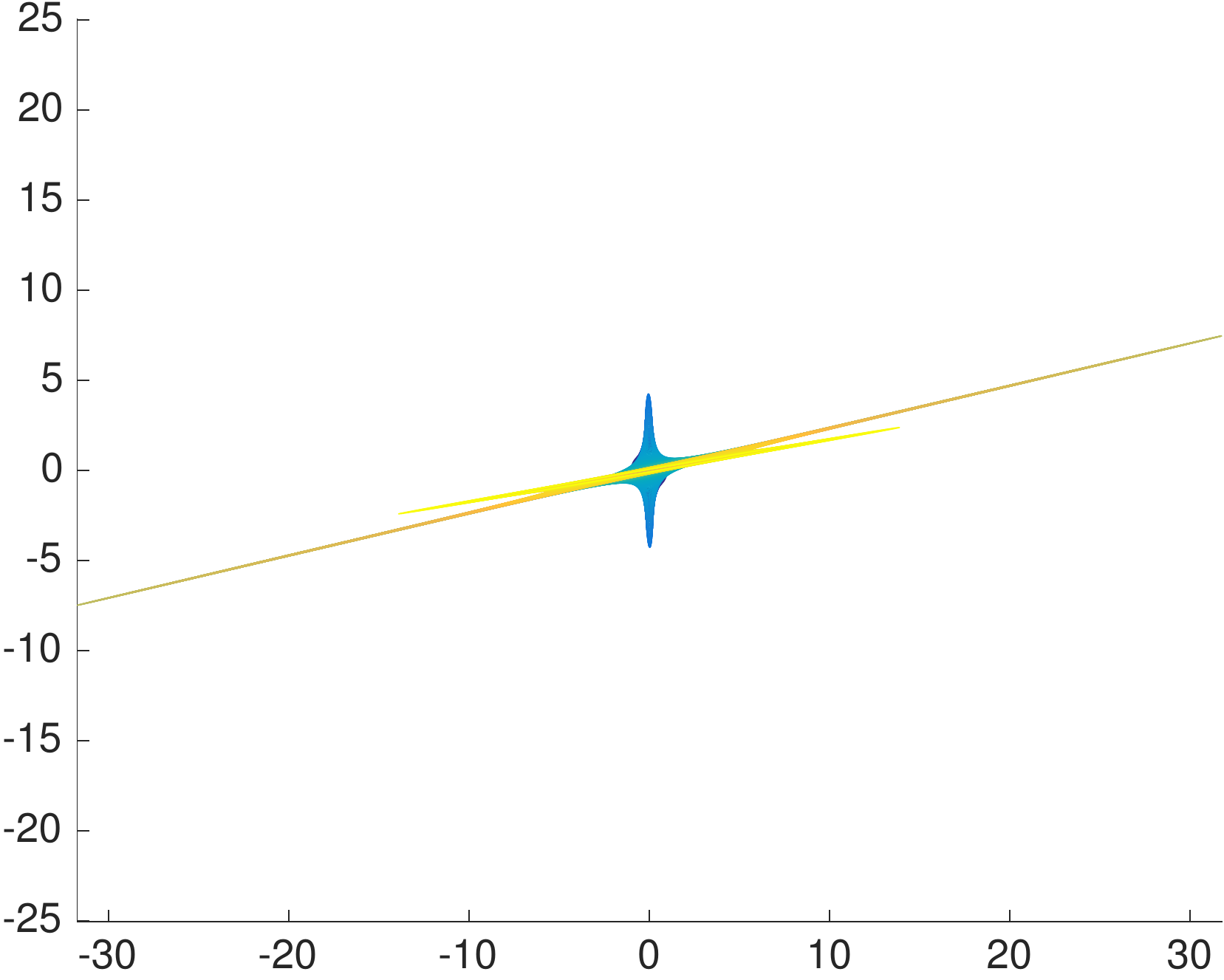}
\hfill
\includegraphics[width = 0.33\textwidth]{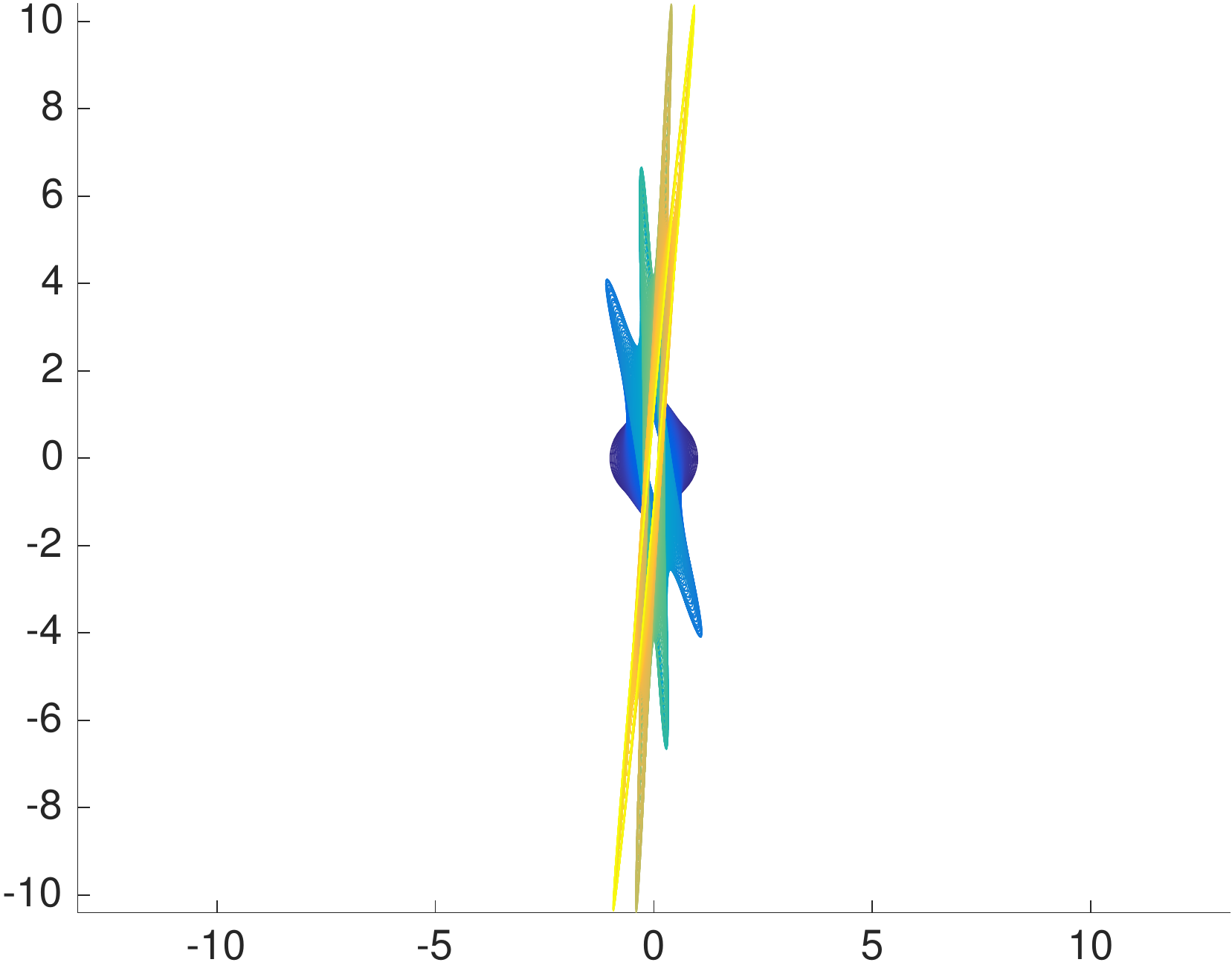}
\caption{The linearized pullback set~$\set{S}(0,T)$ of the periodically driven double gyre system for~$T=5$ and different initial conditions. The color of the contour line corresponds to the time $t$ (blue: initial, yellow: final).}
\label{fig:DGper_galaxies}
\end{figure*}

Figure~\ref{fig:DGper_galaxies} shows~$\set{S}(0,T)$ for $T=5$ and different initial conditions:
\begin{alignat*}{3}
x_{1,0} &= \begin{pmatrix} 1.0 \\ 0.5 \end{pmatrix}, \quad &
x_{2,0} &= \begin{pmatrix} 0.5\\ 0.4 \end{pmatrix}, \quad &
x_{3,0} &= \begin{pmatrix} 0.5\\ 0.7 \end{pmatrix},\\
x_{4,0} &= \begin{pmatrix} 0.86\\ 0.25 \end{pmatrix},\quad &
x_{5,0} &= \begin{pmatrix} 0.99\\ 0.01 \end{pmatrix},\quad &
x_{6,0} &= \begin{pmatrix} 0.98\\ 0.25 \end{pmatrix}.
\end{alignat*}
Recall that the set $\set{S}(0,T)$ is the union over the observation times of pullbacks of the unit circle with respect to the linearized dynamics, and thus we color the single pullbacks with respect to the time they belong to.
The inner white regions are not filled because we only plot the boundaries of the respective ellipses.

\begin{figure*}[htbp]
\centering
\includegraphics[width = 0.33\textwidth]{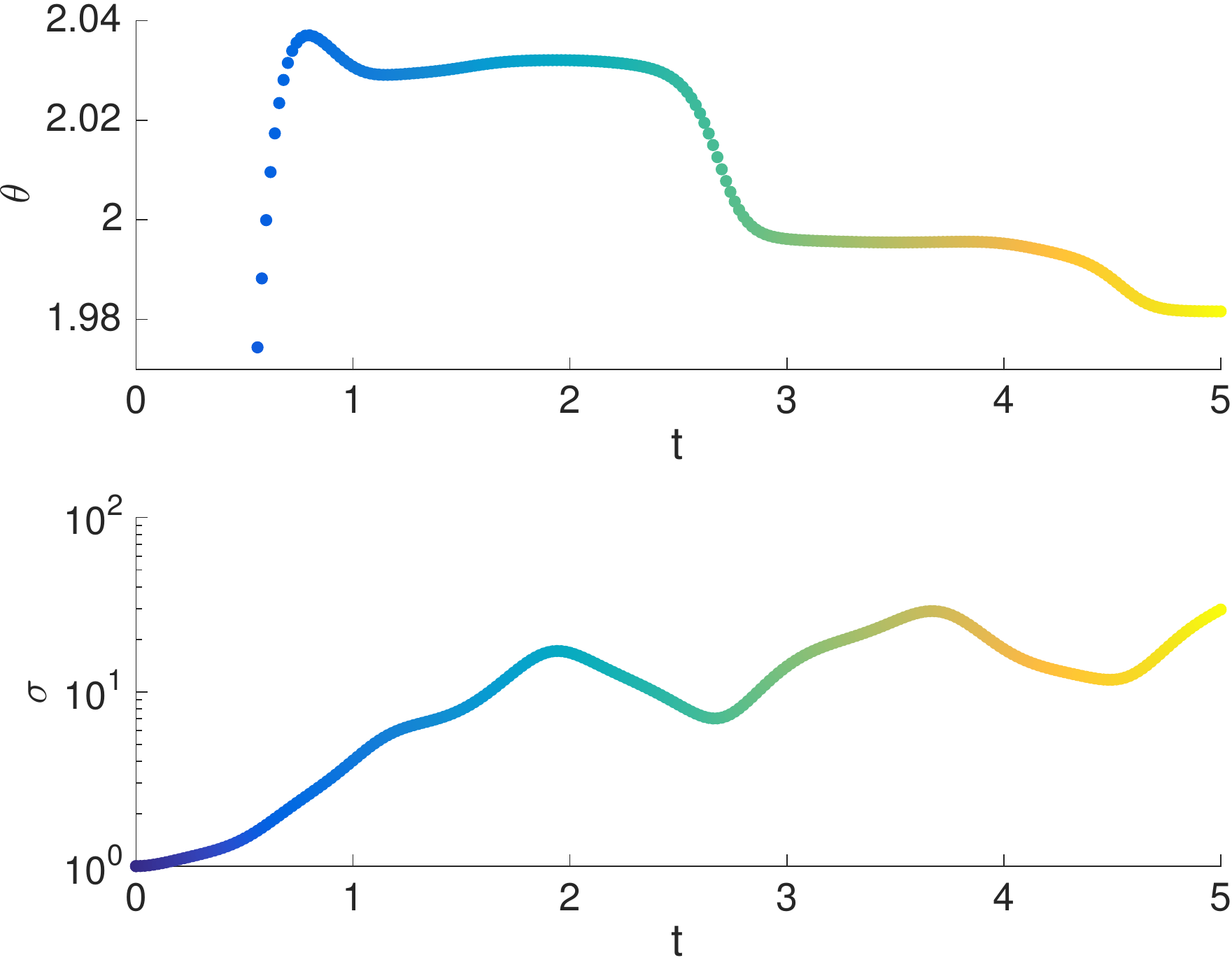}
\hfill
\includegraphics[width = 0.33\textwidth]{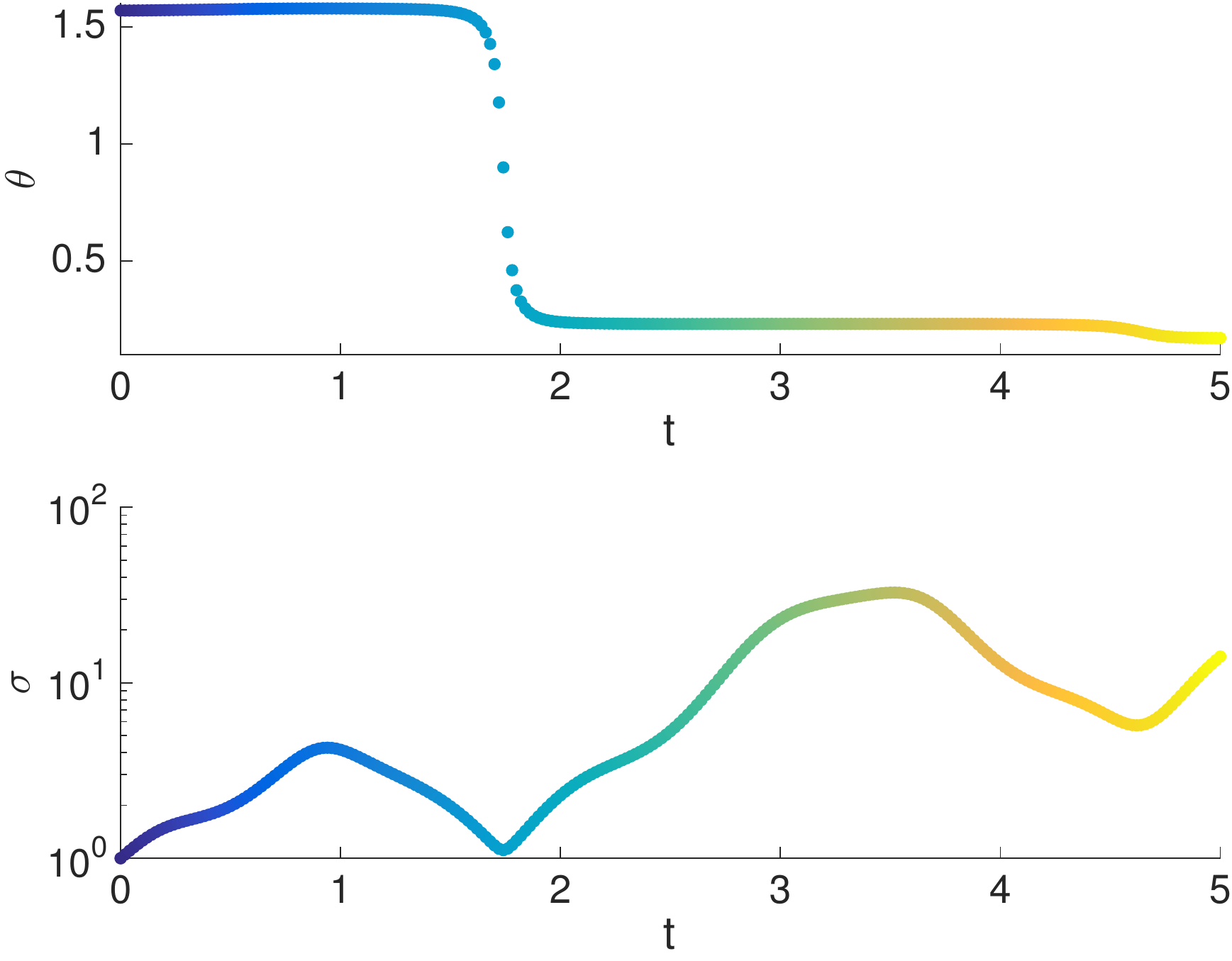}
\hfill
\includegraphics[width = 0.33\textwidth]{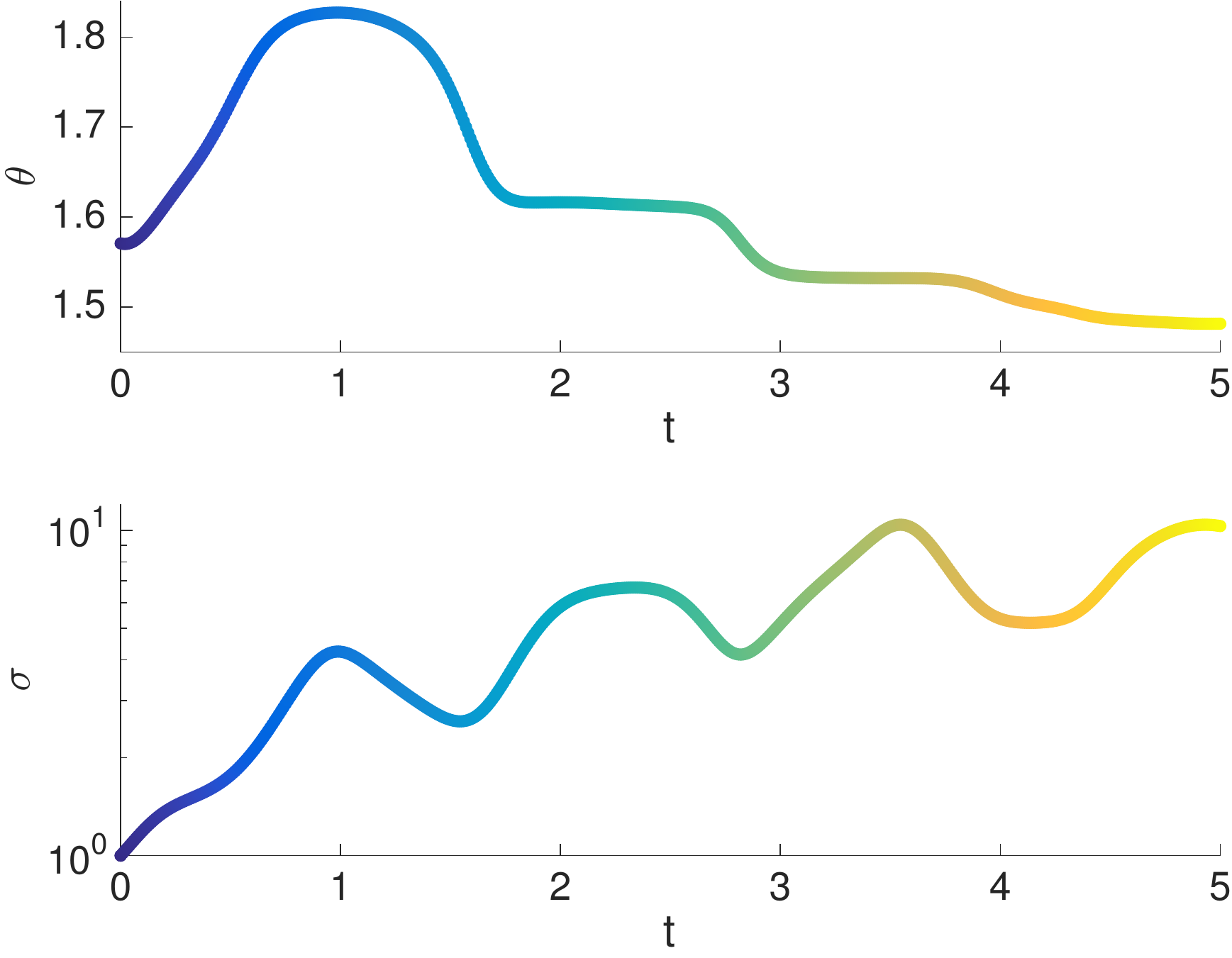}

\caption{Top row: orientation~$\theta(t)$ of the first singular vector of the fundamental matrix~$W(0,t)$ of the double gyre system for~$t\in[0,5]$ for different initial conditions: $x = x_{1,0},x_{5,0}, x_{6,0}$, from left to right. Bottom row: time evolution of the associated singular value~$\sigma(t)$ on a logarithmic scale. The color corresponds to the time $t$ (blue: initial, yellow: final).}
\label{fig:DGper_thetasigma}
\end{figure*}

Assumption~\ref{ass:alignment} extended by the simplified considerations in section~\ref{ssec:degree} seem to be valid: For the trajectories showing strong finite-time hyperbolic behavior ($x_1,x_5$, and $x_6$ moderately) the direction of the pullback ellipses stay constant, or change rapidly in short time intervals when the associated ellipse has small eccentricity (i.e., when~$\sigma_1/\sigma_2 \ngg 1$). This is confirmed by Figure~\ref{fig:DGper_thetasigma}. Therein, the top row shows the time evolution of the orientation of the first singular vector~$v_1(t)$ as an angle in~$[0,2\pi)$ for the trajectories~$x_1, x_5, x_6$ (left to right), while the bottom row shows the time evolution of~$\sigma_1(t)$ on a logarithmic scale. 

\begin{figure*}[htbp]
\centering
\includegraphics[height = 40mm]{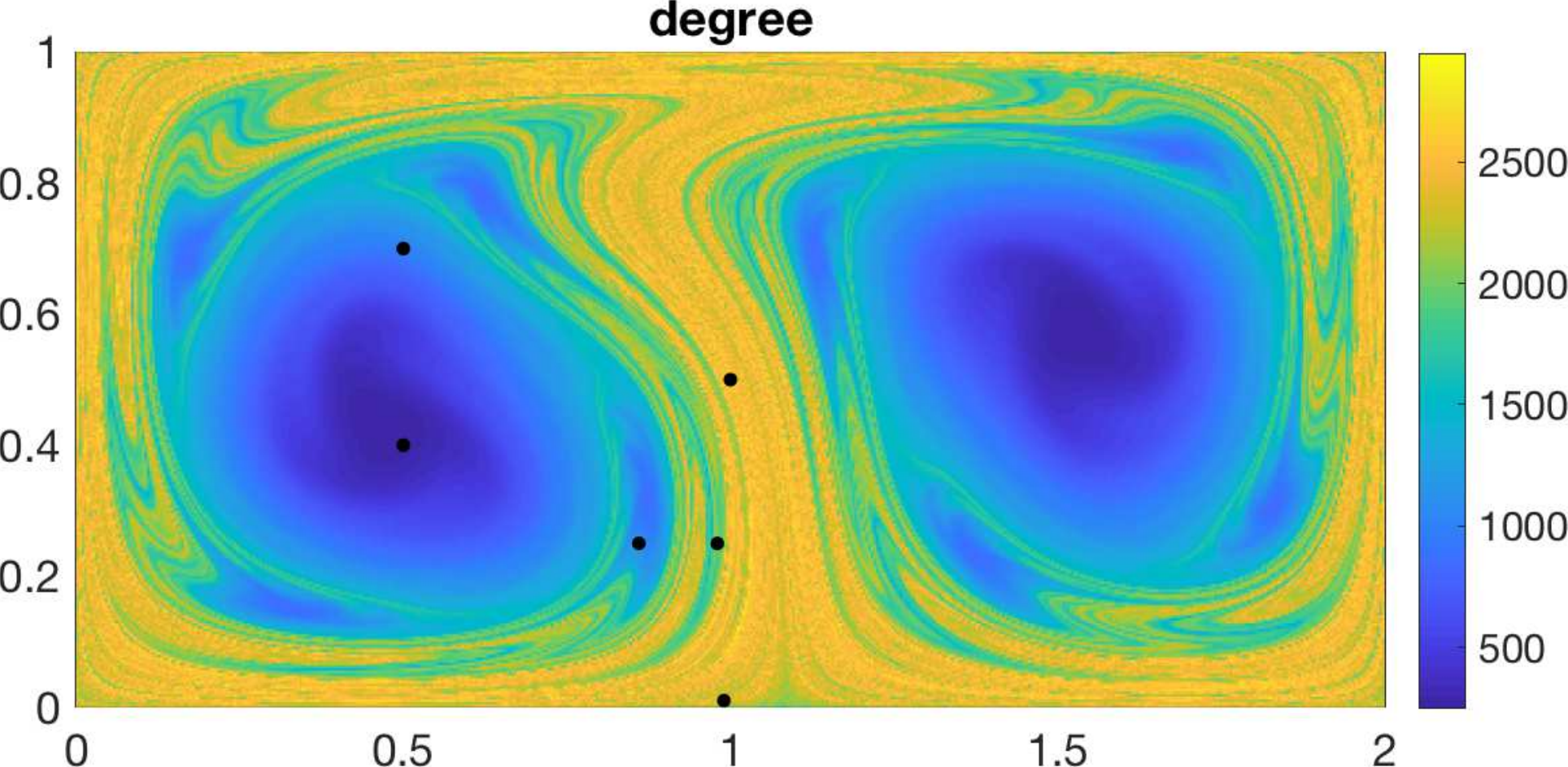}
\hfill
\includegraphics[height = 40mm]{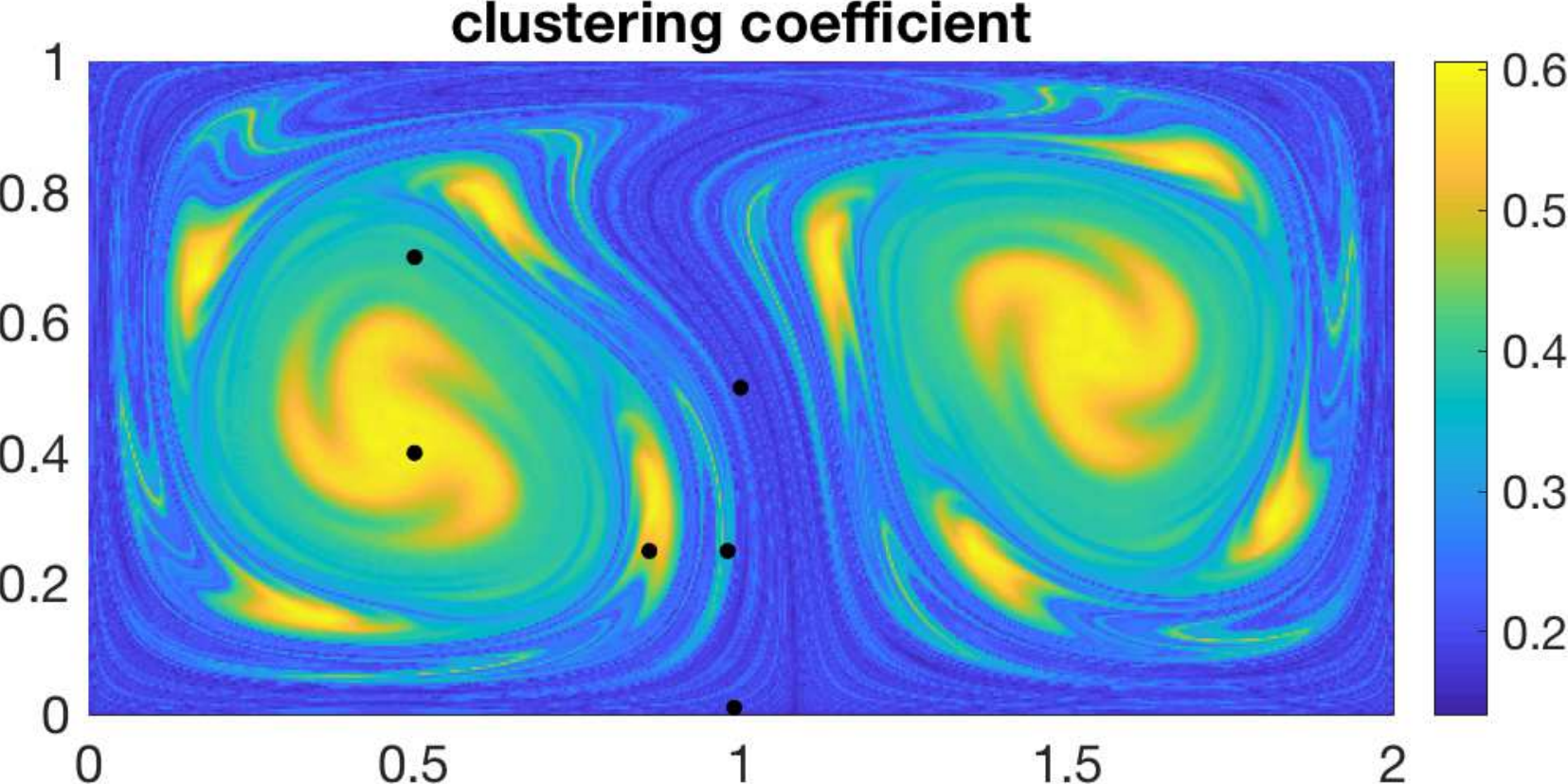}\\
\includegraphics[height = 40mm]{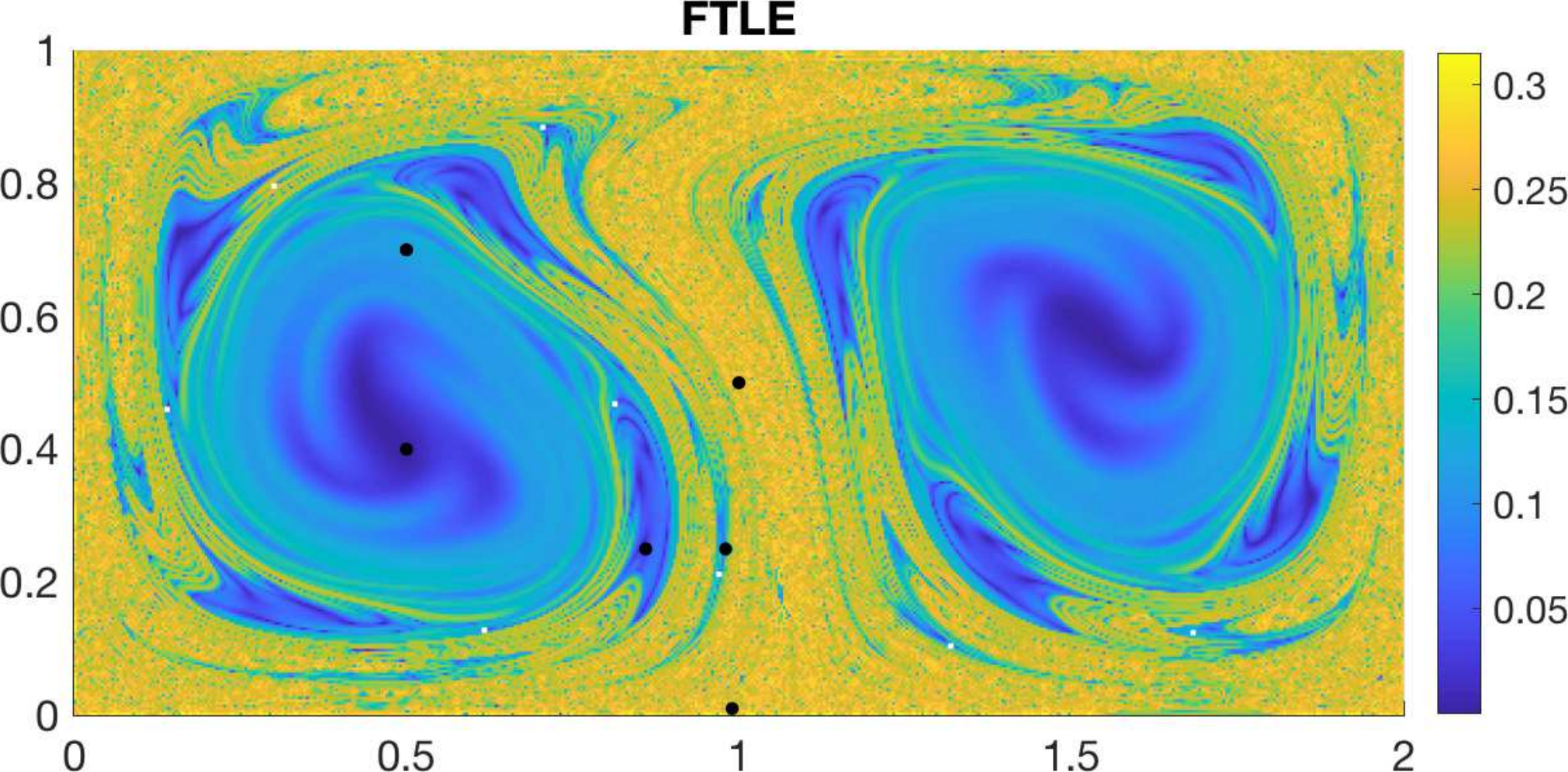}
\hfill
\includegraphics[height = 39mm]{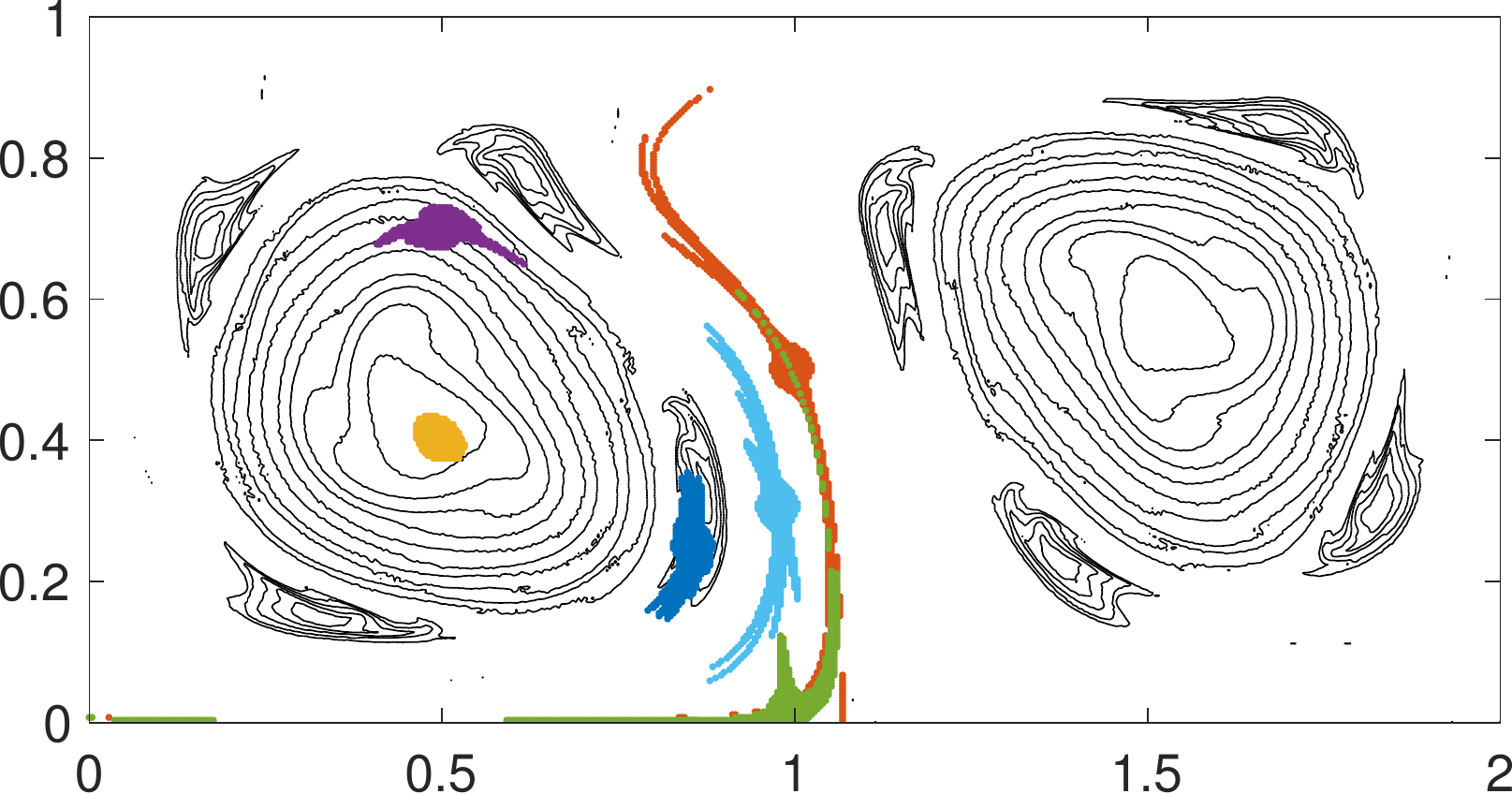}\\
\caption{Degree (top left), clustering coefficient (top right), and FTLE field (bottom left) of the double gyre, with markers for the different initial conditions used in Figure~\ref{fig:DGper_galaxies}. For the degree and clustering coefficient, $\ep = 0.03$ and flow time $T=20$ with $\Delta t = 0.1$ was used. Bottom right: The initial conditions of all trajectories that come $\ep$-close to one of the~$x_i$, $i=1,\ldots,6$ for $T=5$, i.e., $\set{G}_{\ep}(x_i)$ (the shorter time span was taken for better visual comparability). }
\label{fig:DG_degcc}
\end{figure*}

We also construct a network with $500\times 251$ trajectories spaced equally on a grid in the domain~$[0,2]\times [0,1]$. 
We choose $\ep = 0.03$, flow time $T = 20$ and $\Delta t = 0.1$. Degree, clustering coefficient, and the sets $\smash{ \bigcup_{t\in[0,5]}  \phi(t_0,t)^{-1} \set{B}_{\ep}(x(t)) }$ for~$x_1,\ldots,x_6$, of which the~$\set{S}(0,5)$ are linearized approximations, are shown in Figure~\ref{fig:DG_degcc}. 

\begin{figure}[htbp]
\includegraphics[width = 0.35\textwidth]{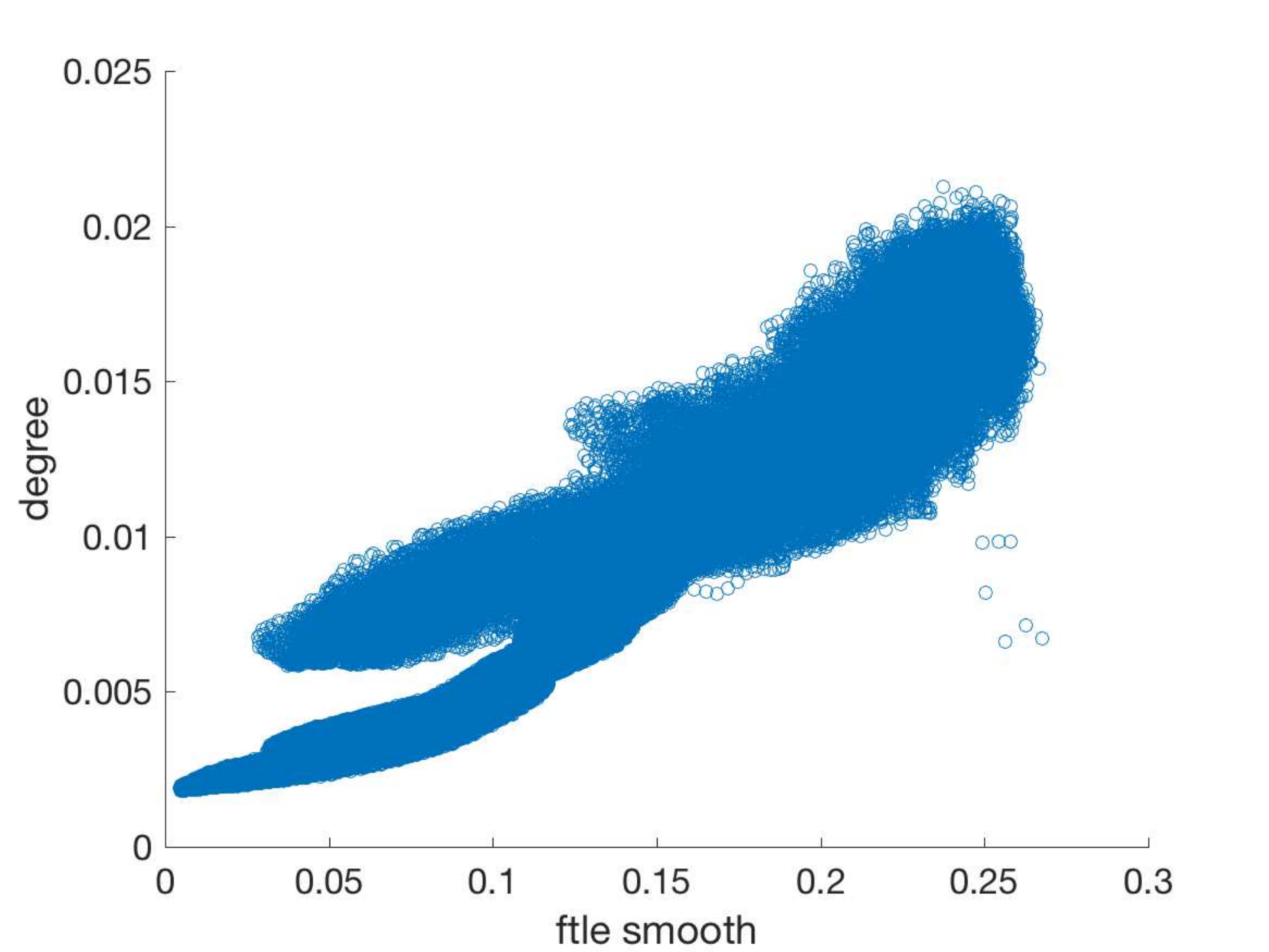}\hfill
\includegraphics[width = 0.35\textwidth]{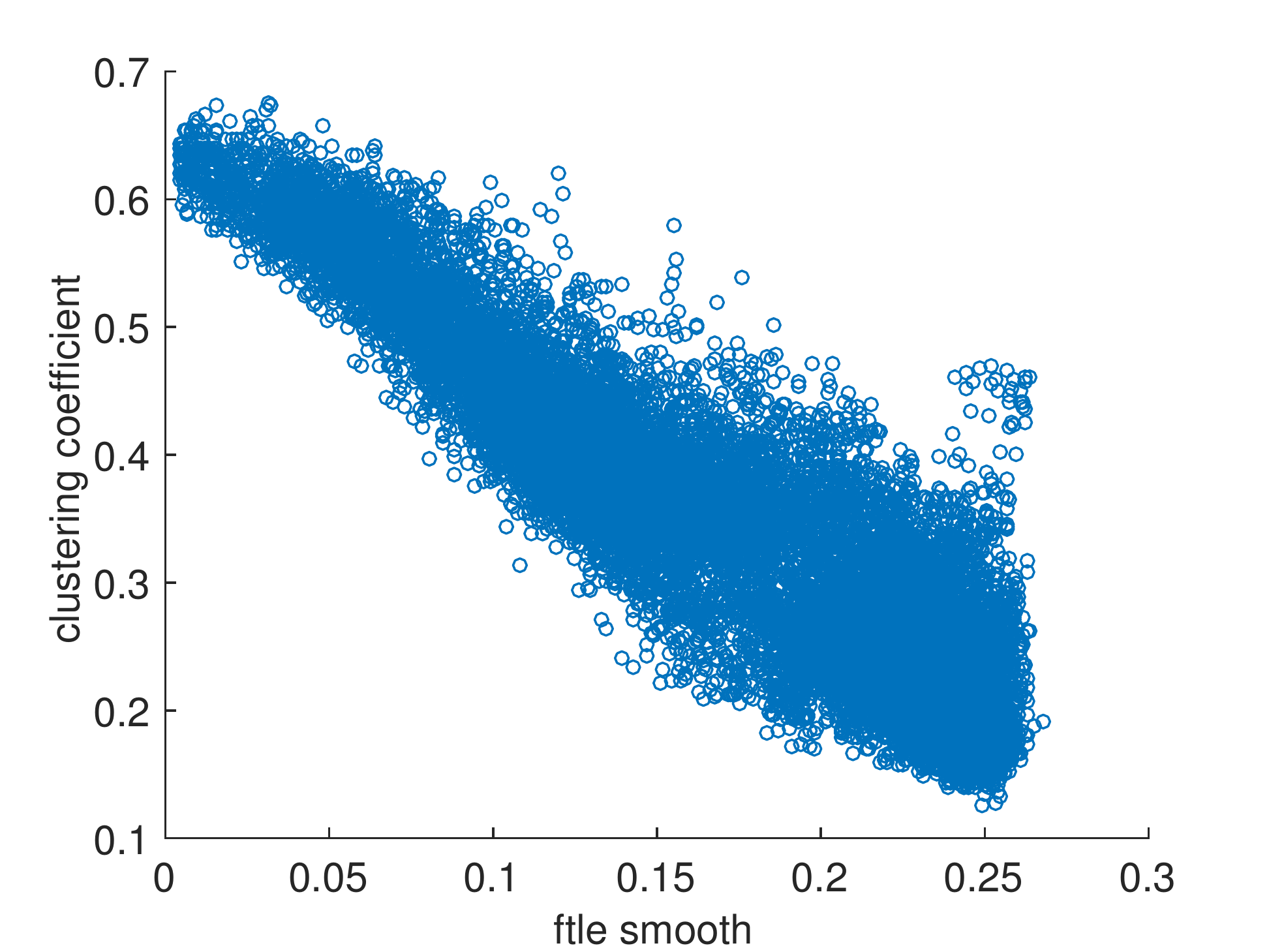}
\caption{Top: smoothed FTLE field vs.\ the relative degree (degree normalized by~$n$). The correlation coefficient is $0.95$. Bottom: smoothed FTLE field vs. the local clustering coefficient. The correlation coefficient is $-0.93$. }
\label{fig:DG_ftle}
\end{figure}

For comparison we also compute the FTLE field \eqref{eq:ftle} with the same resolution, see Figure \ref{fig:DG_degcc} (bottom left). We see that the affine-linear relation between FTLE and degree, predicted by Proposition~\ref{prop:vol_galaxy}  in an idealized (linearized) setting, holds only up to a substantial spread in the values; see Figure~\ref{fig:DG_ftle}. Some quantitative agreement is clearly visible though, and the correlation coefficient between them is $0.95$ for a smoothed FTLE-field.\footnote{Structures that are on scales of higher order in $\ep$ are not captured by the degree anyway, and the linearized relations~\eqref{eq:pullbackdist} and~\eqref{eq:pullbackballs} disregard errors on higher than first order too; without smoothing the FTLE field, it is~$0.88$.}
The clustering coefficient and FTLE are negatively correlated with a correlation coefficient of $-0.93$, indicating that the region with chaotic dynamics is already in a filamenting regime, but did not yet reach well-mixedness; cf.\ Table~\ref{tab:regimes}. In this example we do not consider the measure ``closeness'', as it does not contribute additional insights.

\begin{figure*}[htbp]
\includegraphics[width = 0.49\textwidth]{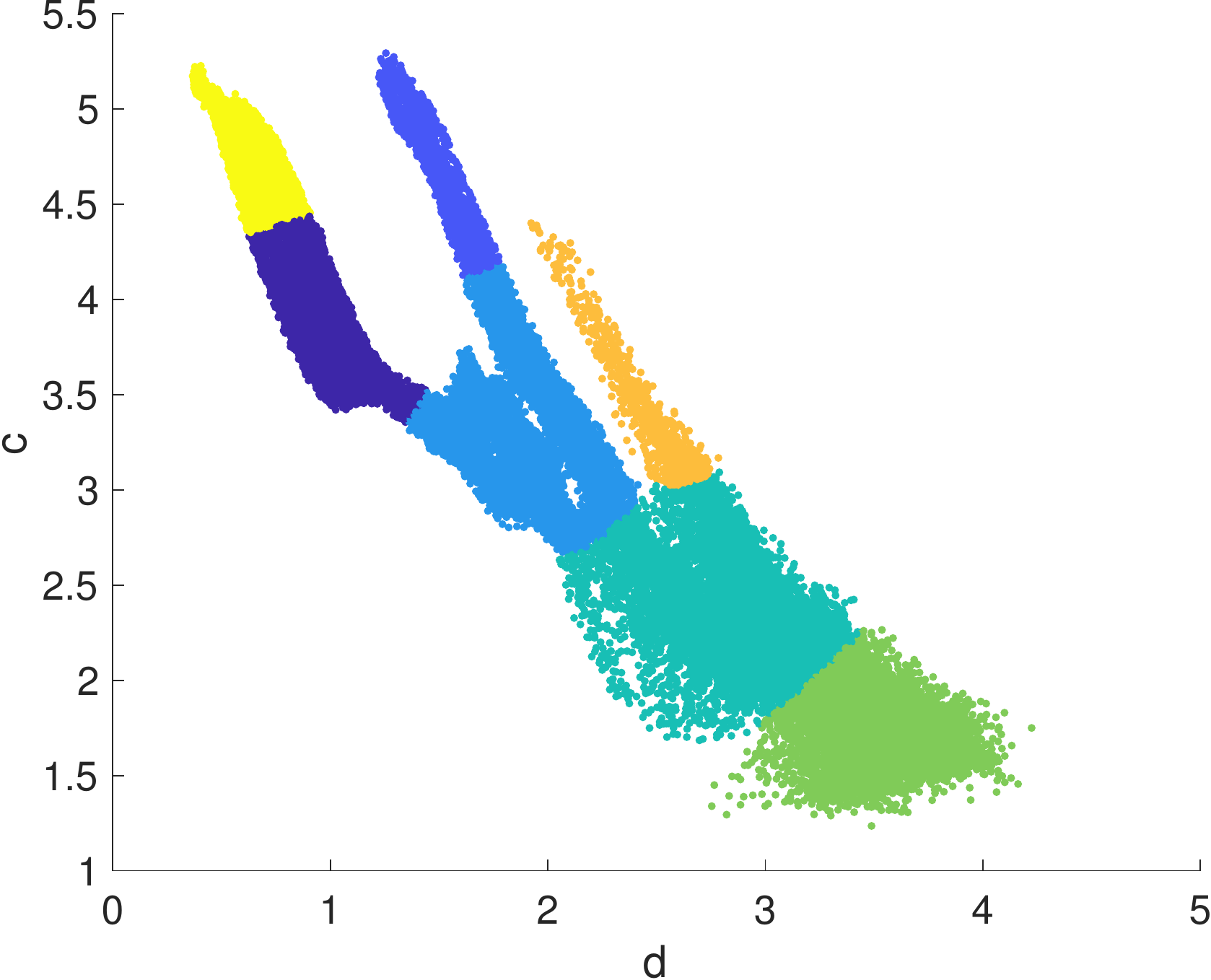}\hfill
\includegraphics[width = 0.49\textwidth]{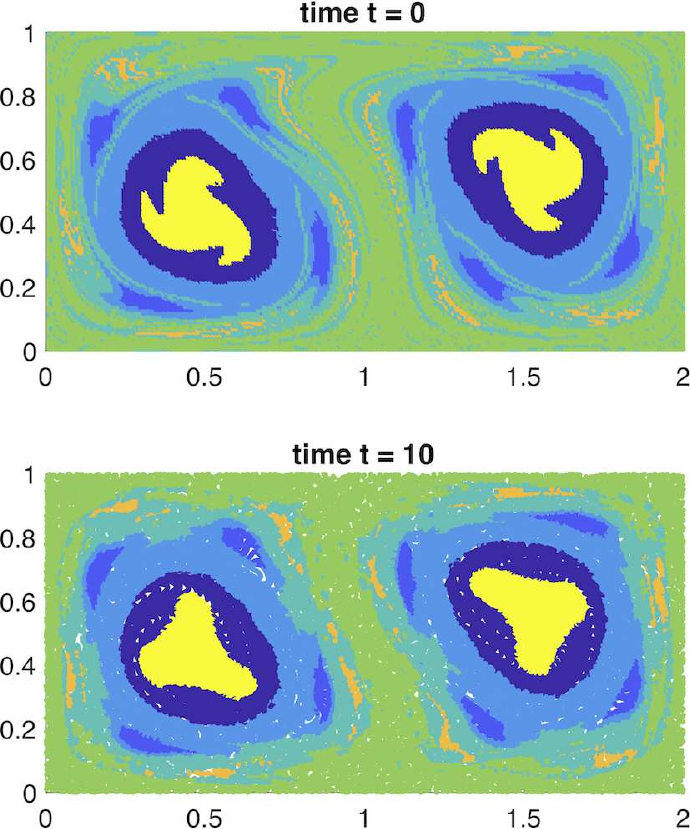}
\caption{Left: When the degree $d$ of trajectories is plotted against the clustering coefficient $c$, one can see that they are not perfectly correlated: The points with large $c$ fall into three branches. The colors are a classification produced with diffusion maps and subsequent clustering by $k$-means\cite{k-means}. Right: Same classes as on the left, trajectories plotted in state space at initial time~$t=0$ and halftime $t=10$.}
\label{fig:DG_classes}
\end{figure*}

Next we investigate which dynamical structures can be identified and distinguished from another by looking at the two network measures, degree and clustering coefficient, \emph{simultaneously}. To this end we consider the point cloud
\begin{equation} \label{eq:embedding2d}
\set{E} = \left\{ (\bar{d}_i, \bar{c}_i)\,\big\vert\, i=1,\ldots,n \right\} \subset \R^2,
\end{equation}
where $\bar{d}$ and $\bar{c}$ are the degree and clustering coefficient normalized by their respective standard deviations such that they cover a comparable numerical range,
and analyze this set by the established manifold-learning tool, the so-called \emph{diffusion maps}~\cite{CoLa06} together with clustering. Diffusion maps, in a nutshell, finds intrinsic coordinates on a point-cloud approximation of a manifold, such that these coordinates are monotonic in the geodesic distance along the manifold. As an effect, if the point cloud has a complicated topology in its original space, the diffusion-map coordinates tend to ``disentangle'' it, and clustering in this new space reveals regions of the point cloud that are close-by with respect to the \emph{intrinsic geodesic} distance of this set. This is shown in Figure~\ref{fig:DG_classes} for proximity parameter $\epsilon=0.01$ in the diffusion maps algorithm and clustering its seven dominant eigenvectors into seven clusters.

\begin{figure}[htbp]
\centering
\includegraphics[width = 0.49\textwidth]{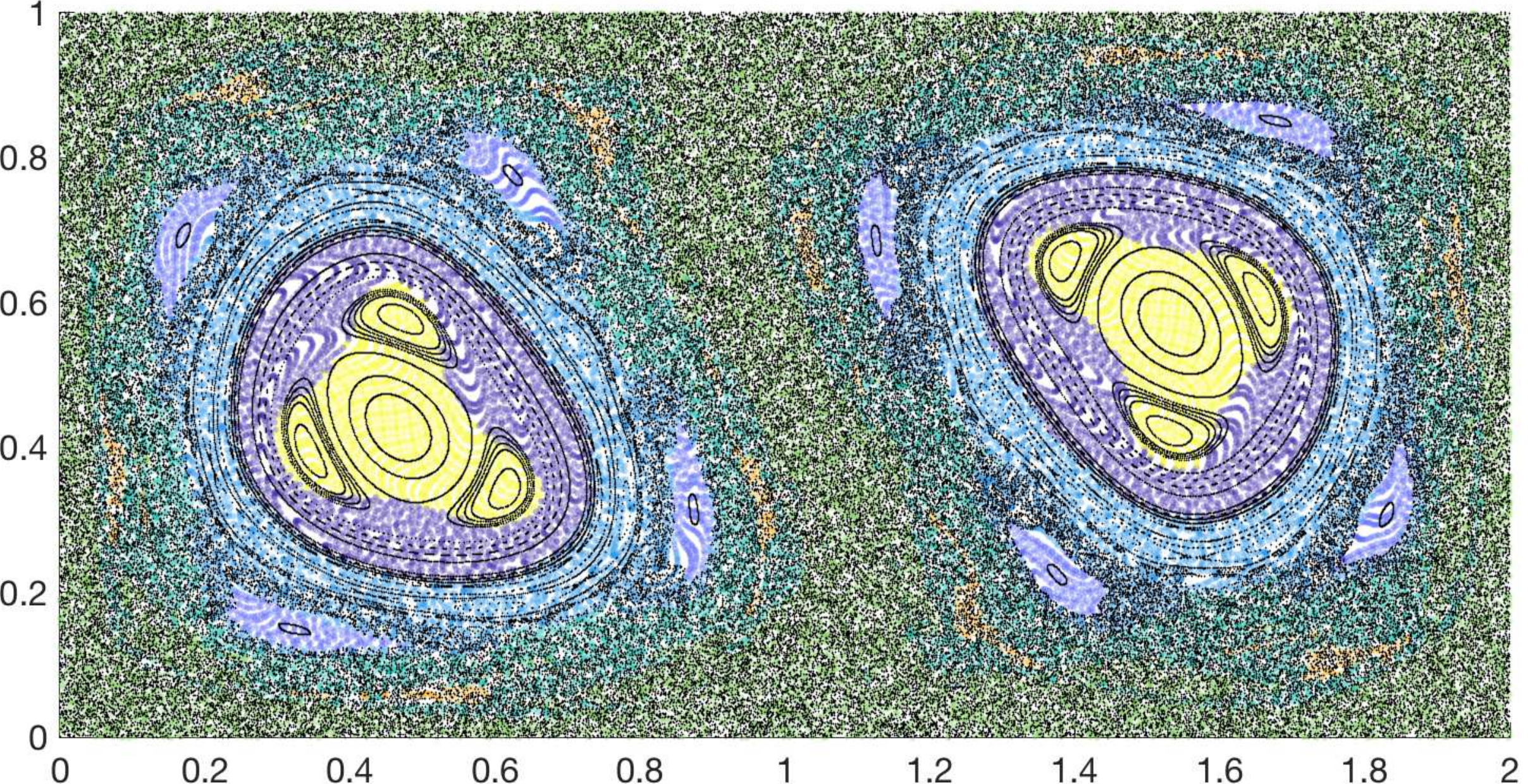}
\caption{
The classification from Figure~ \ref{fig:DG_classes} with trajectory positions at halftime $t=10$ superimposed with a Poincar\'e plot. The classification identifies gyre cores, KAM tori and the chaotic region.}
\label{fig:DG_poincare}
\end{figure}

It is interesting, that the classification of trajectories by their proximity in the set~$\set{E}$ (with respect to its geodesic distance), can be connected to different qualitative dynamical behavior. This is shown in Figure~\ref{fig:DG_poincare}, where the classification is compared with a ``Poincar\'e plot'' of the double gyre system. The classification separates regular regions in the gyre core, KAM tori, and the chaotic region around them.
Also, the invariant ``inner'' and ``outer'' gyre cores are distinguished. Although they all consist of trajectories evolving on invariant cycles, the outer cycles are longer. On the one hand, since trajectories on close-by but different cycles do not keep in phase, the longer cycles have a larger $\ep$-neighborhood, thus a larger degree. On the other hand, on intermediate time intervals (like ours, with $T=20$) this means a smaller clustering coefficient.

\begin{figure}[htbp]
\centering

\includegraphics[width = 0.3\textwidth]{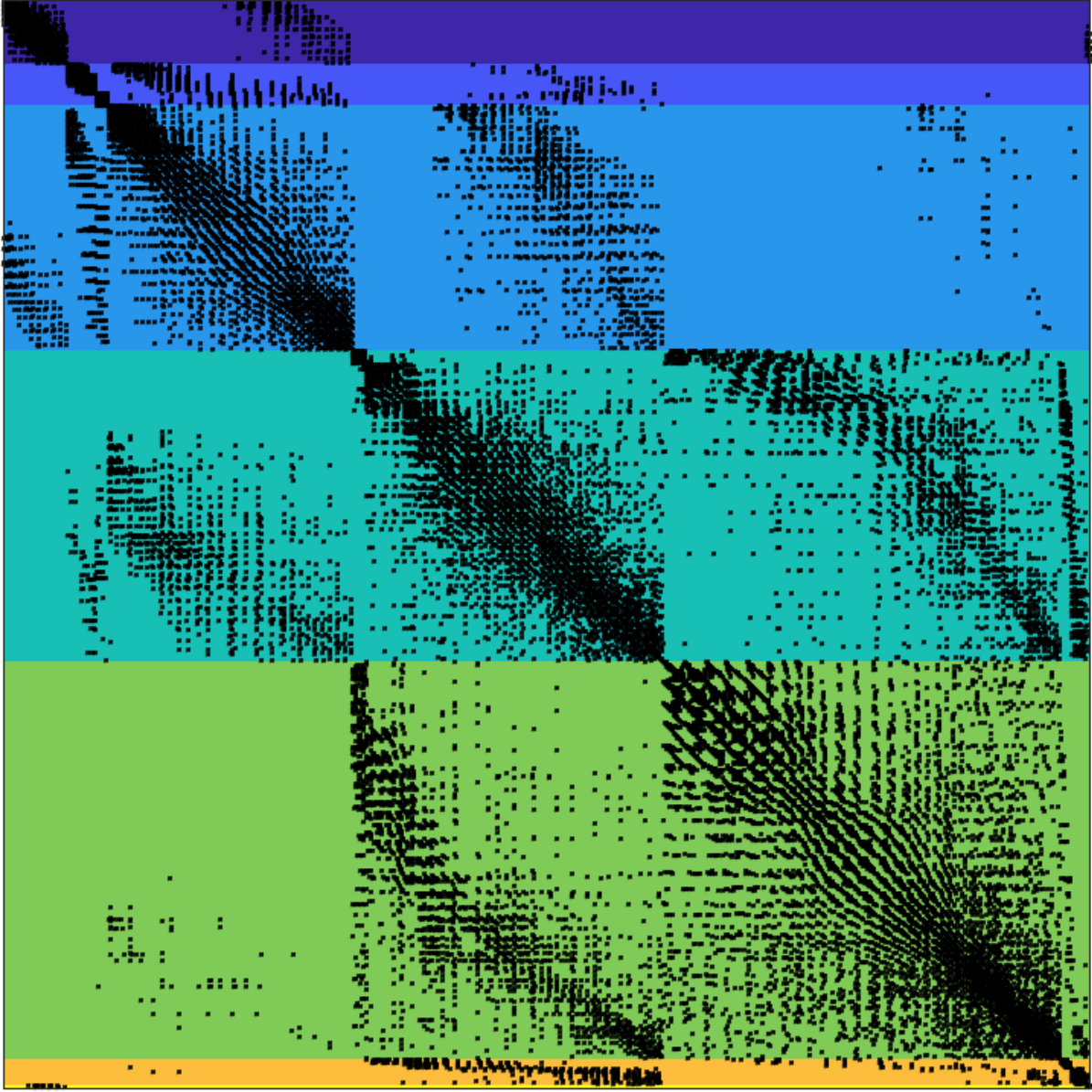}
\hfill
\includegraphics[width = 0.3\textwidth]{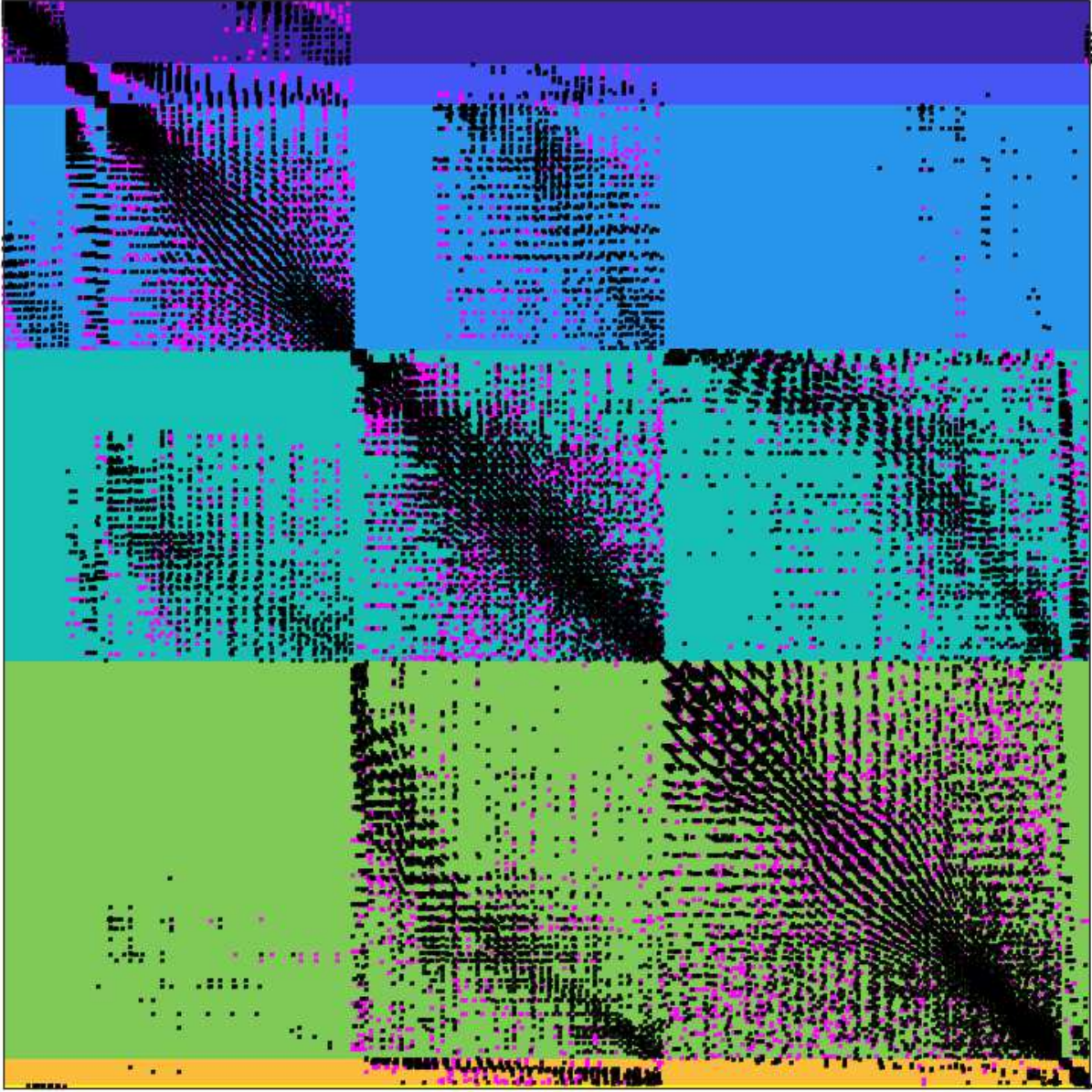}
\hfill
\includegraphics[width = 0.3\textwidth]{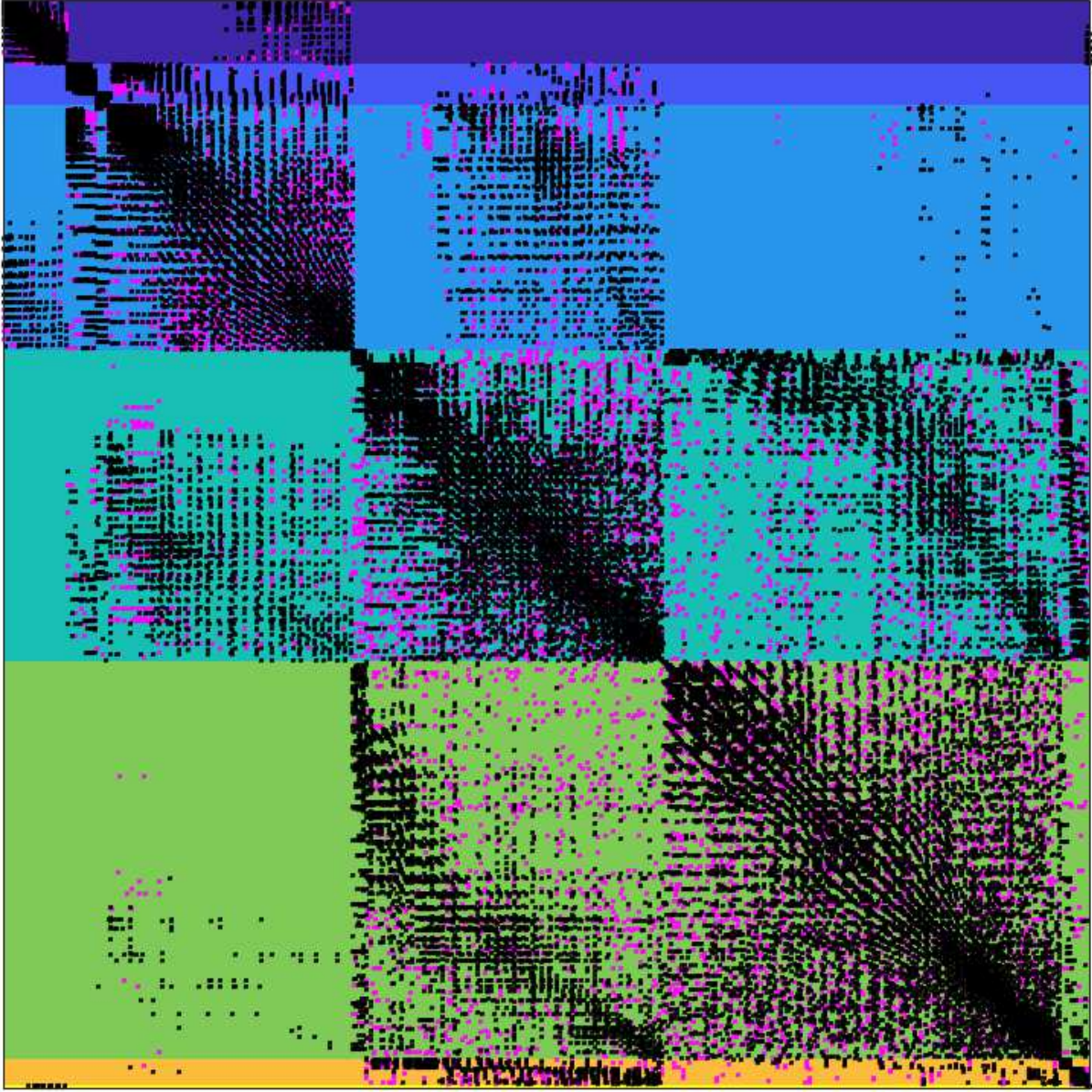}
\caption{
Occupation structure (non-zero entries are show by magenta and black dots) of a $1000\times 1000$ submatrix of the adjacency matrix $A$ for final times $T=5,10,20$, from left to right. The background color of a row is according to the clustering of $\set{E}$, and the same as in the previous figures. The entries that are new compared with the previous final-time matrix are shown in magenta. We note that the largest growth of neighbors for increasing final time $T$ is experienced by the nodes in the regions mixing most strongly, i.e., green and teal.
}
\label{fig:DG_adjacency}
\end{figure}

To gain some additional intuition of the structure of the network and its temporal change, we visualize a part of the adjacency matrix $A$ for the end times~$T=5,10,20$. We do not show the entire matrix due to its size. Please refer to Figure~\ref{fig:DG_adjacency}.

\subsection{Ocean flow}
\label{ssec:aviso}

\begin{figure*}[htbp]
\centering
\includegraphics[width = 0.24\textwidth]{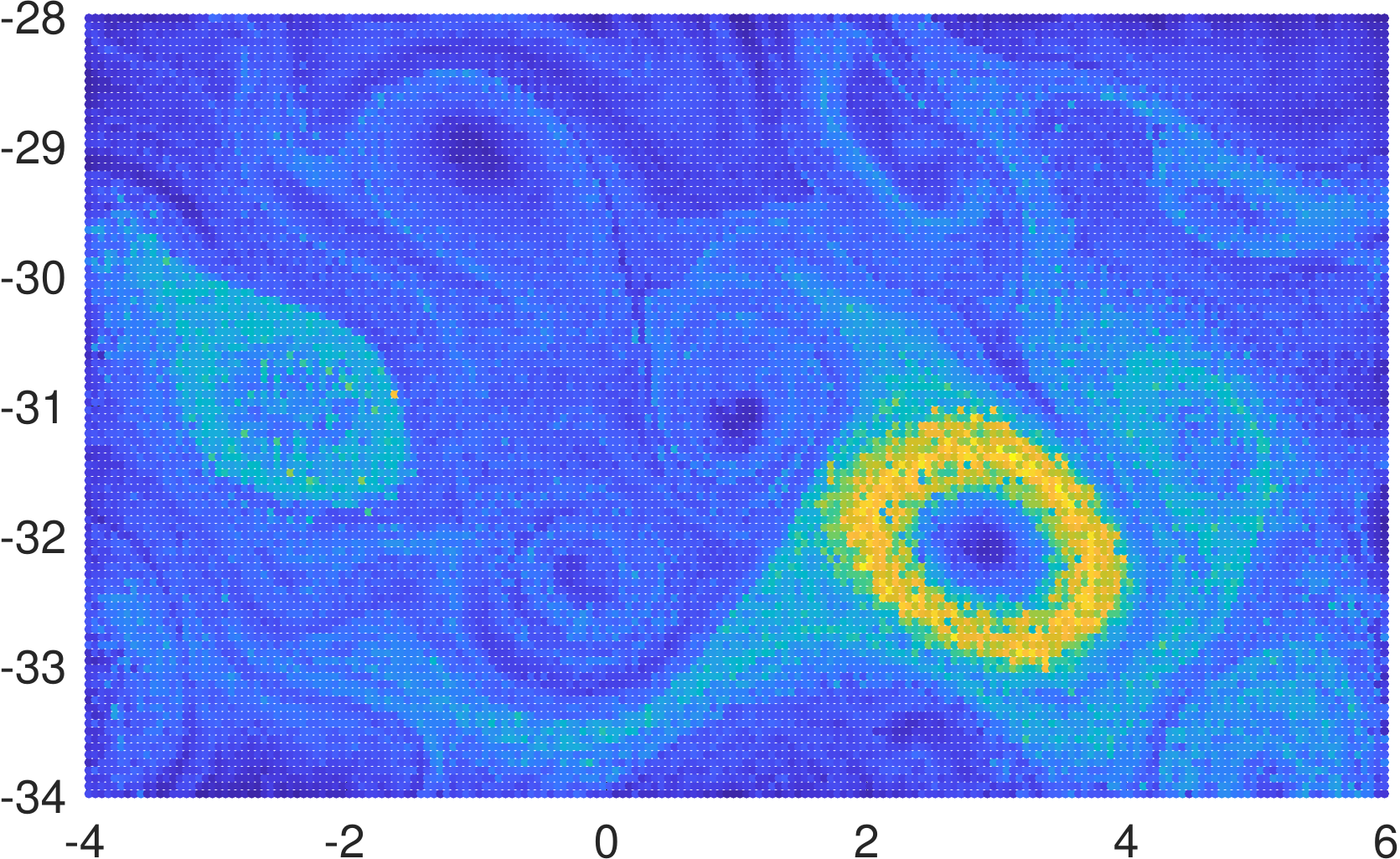}
\hfill
\includegraphics[width = 0.24\textwidth]{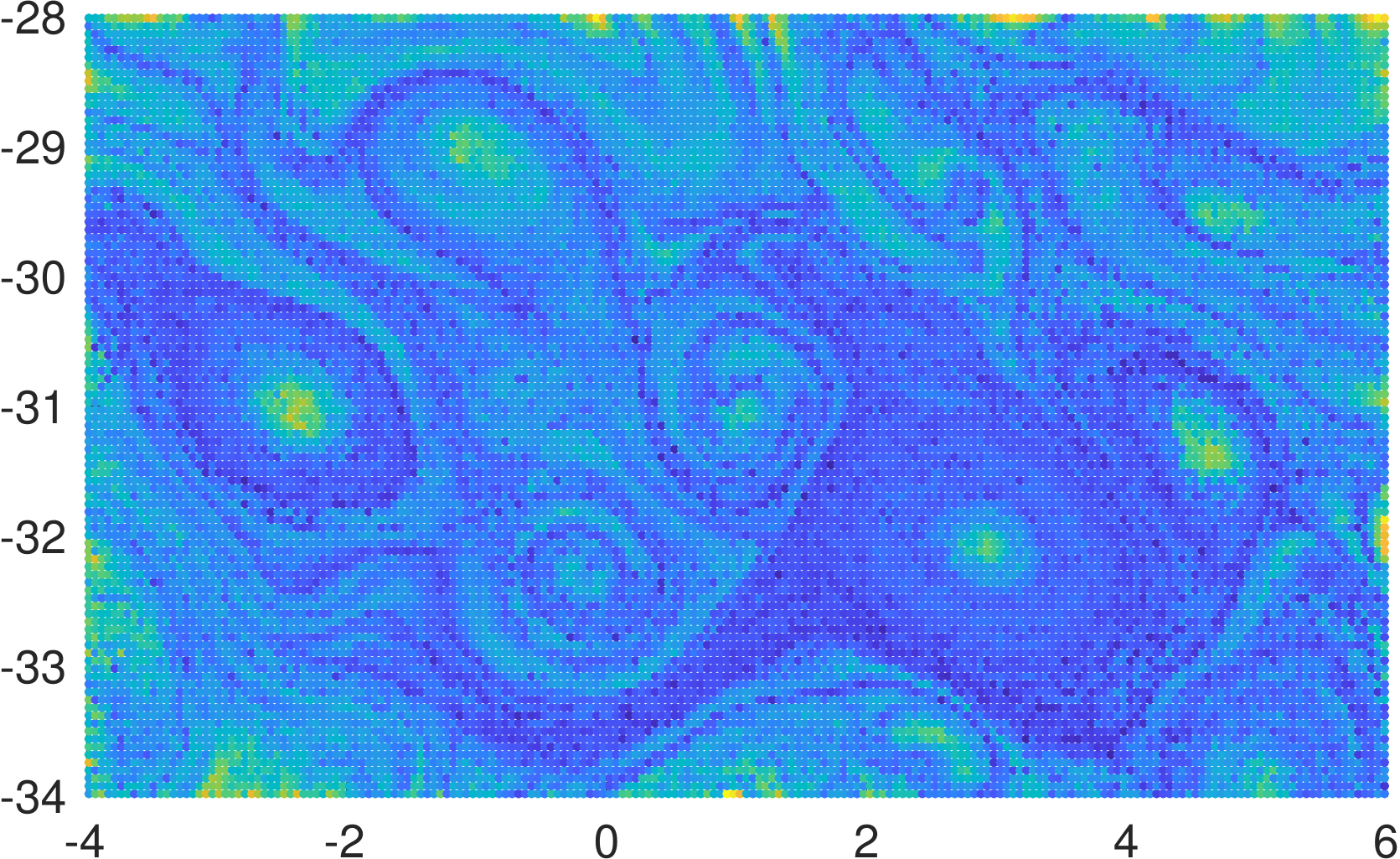}
\hfill
\includegraphics[width = 0.24\textwidth]{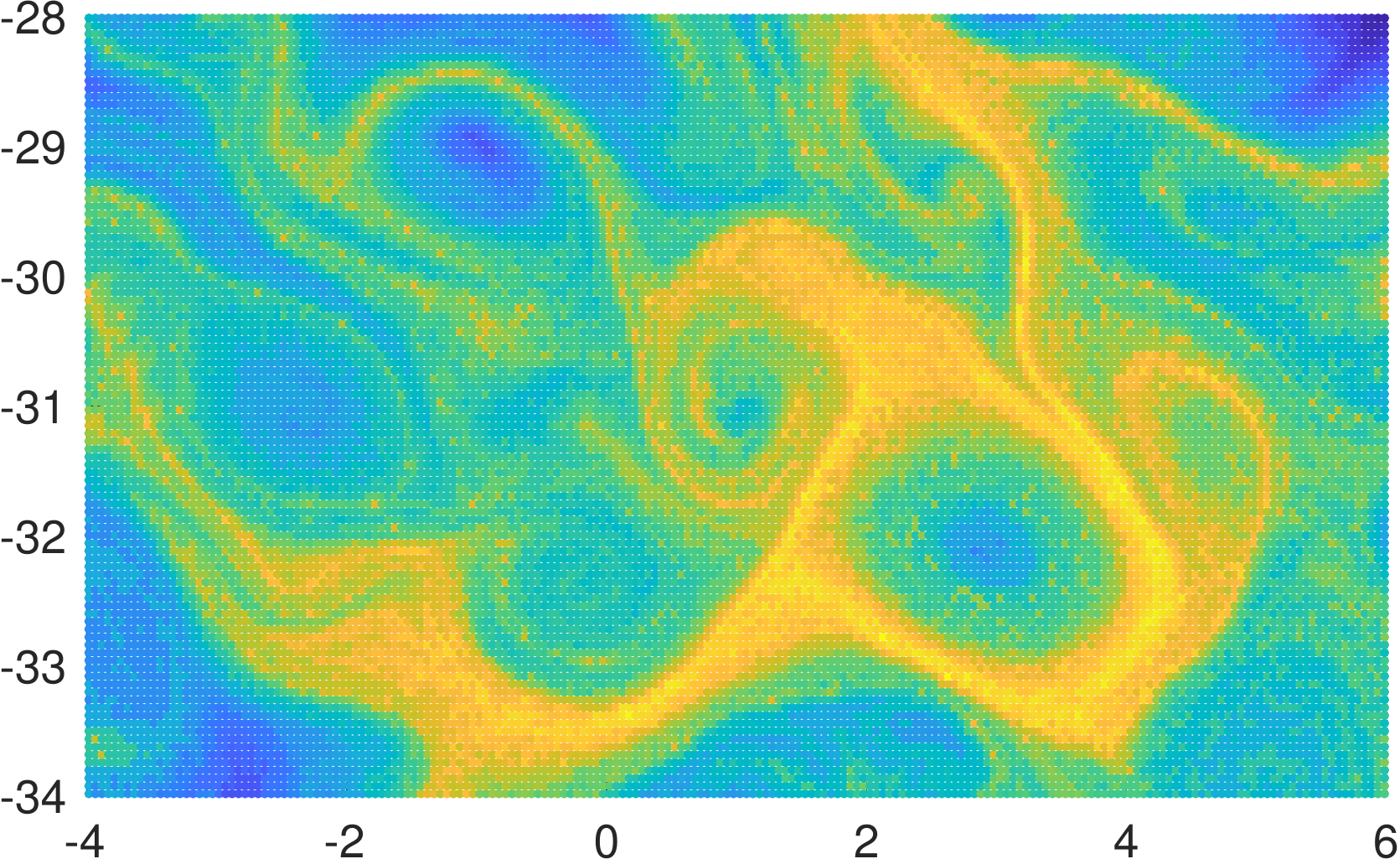}
\hfill
\includegraphics[width = 0.24\textwidth]{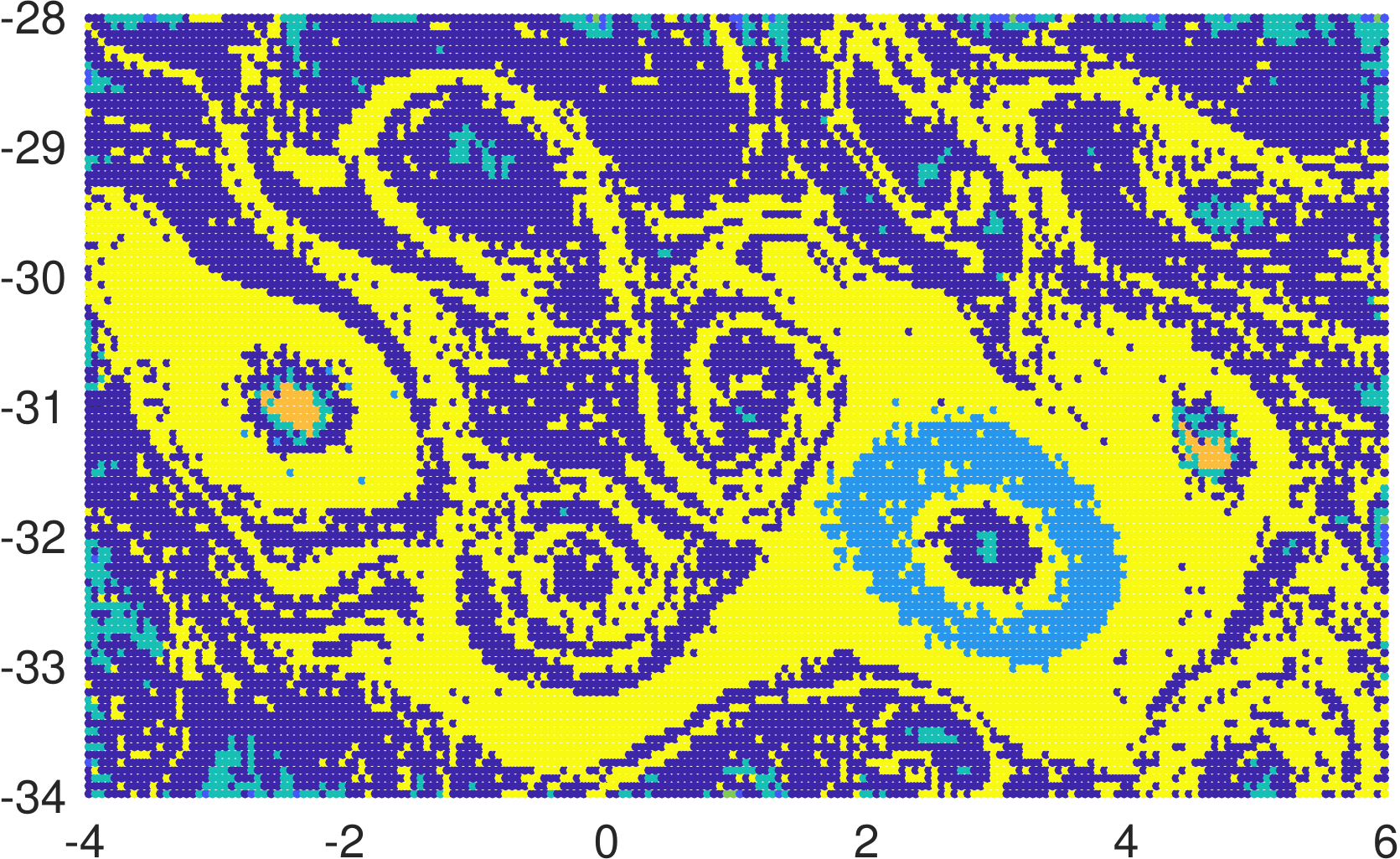}

\includegraphics[width = 0.24\textwidth]{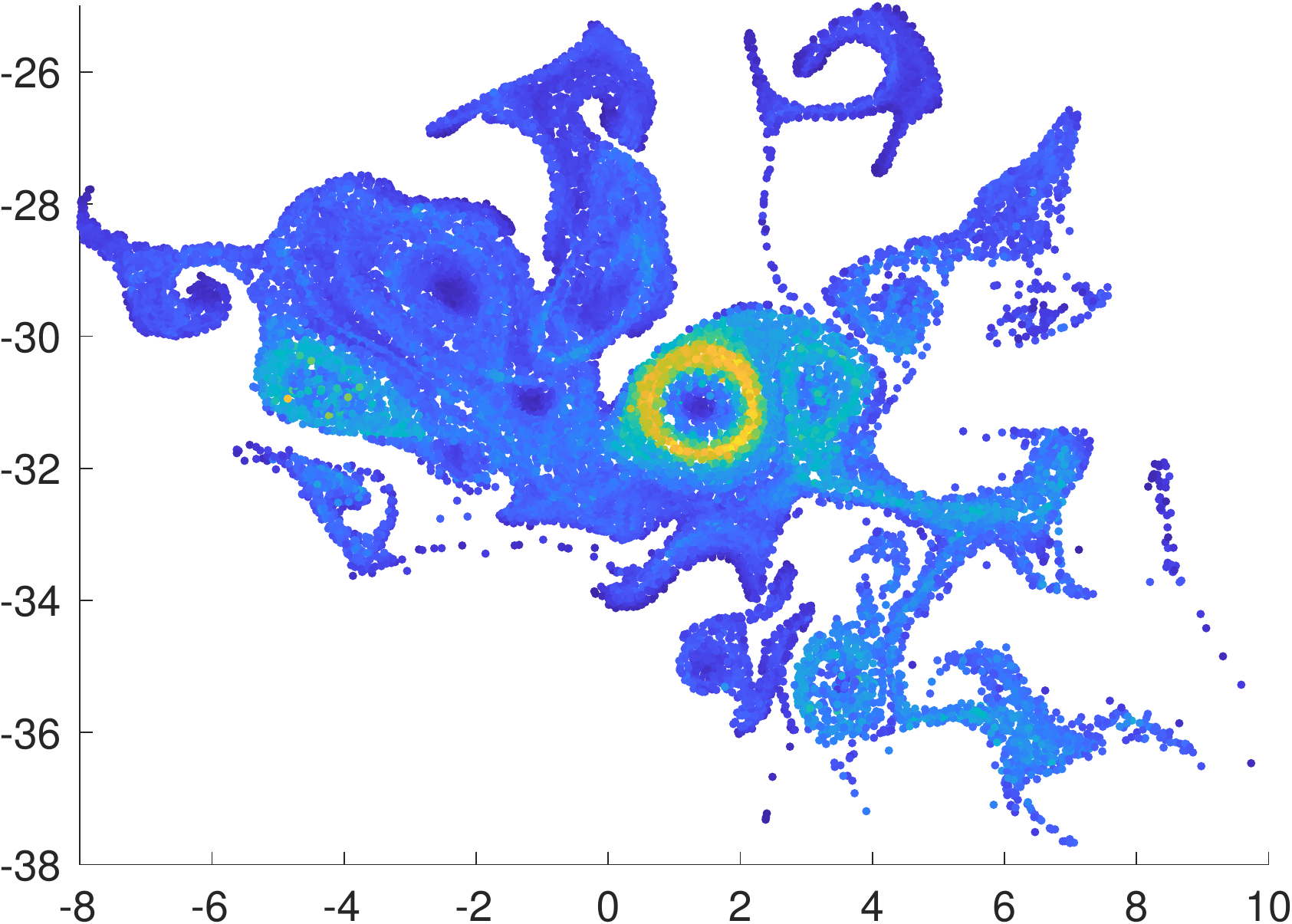}
\hfill
\includegraphics[width = 0.24\textwidth]{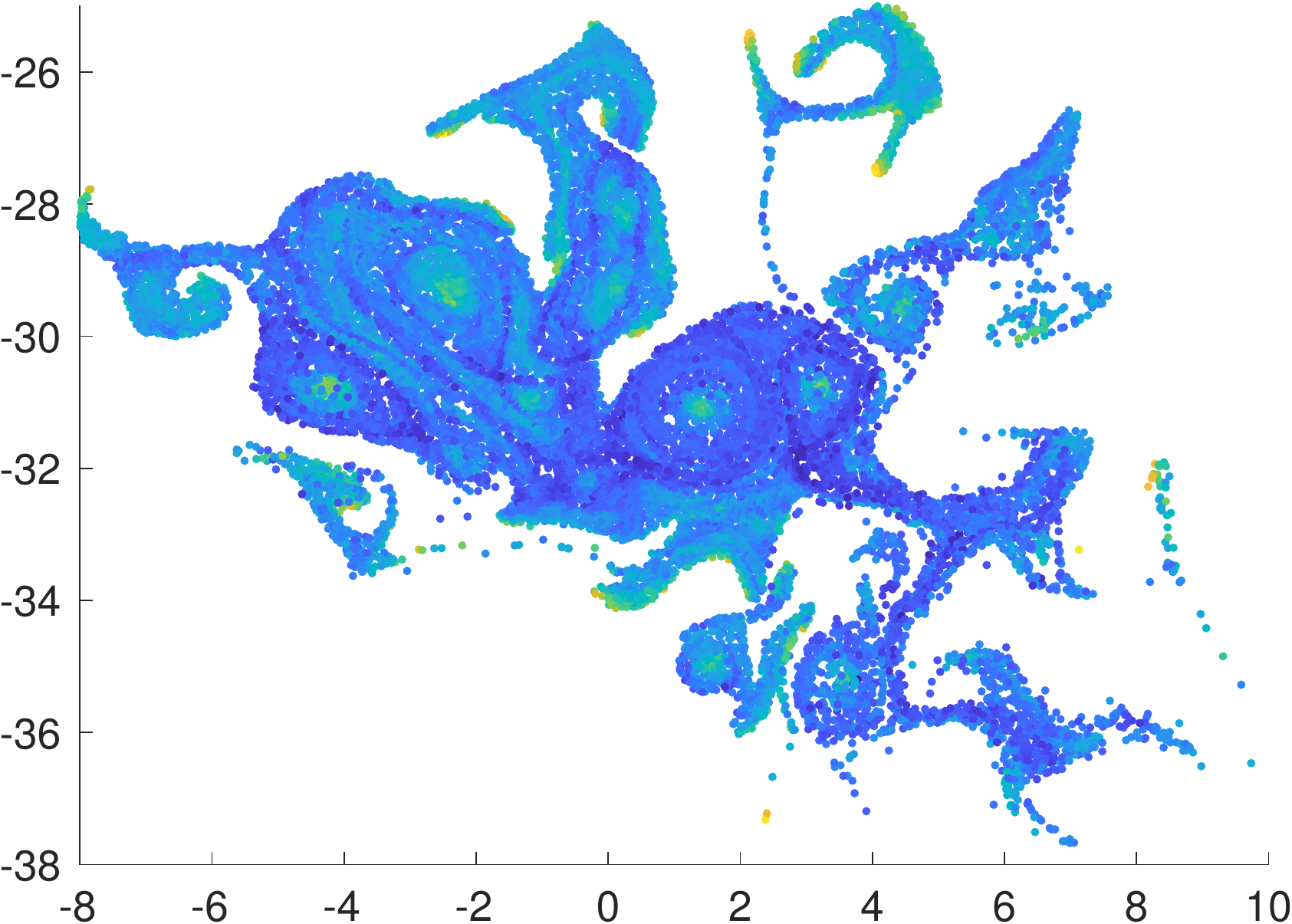}
\hfill
\includegraphics[width = 0.24\textwidth]{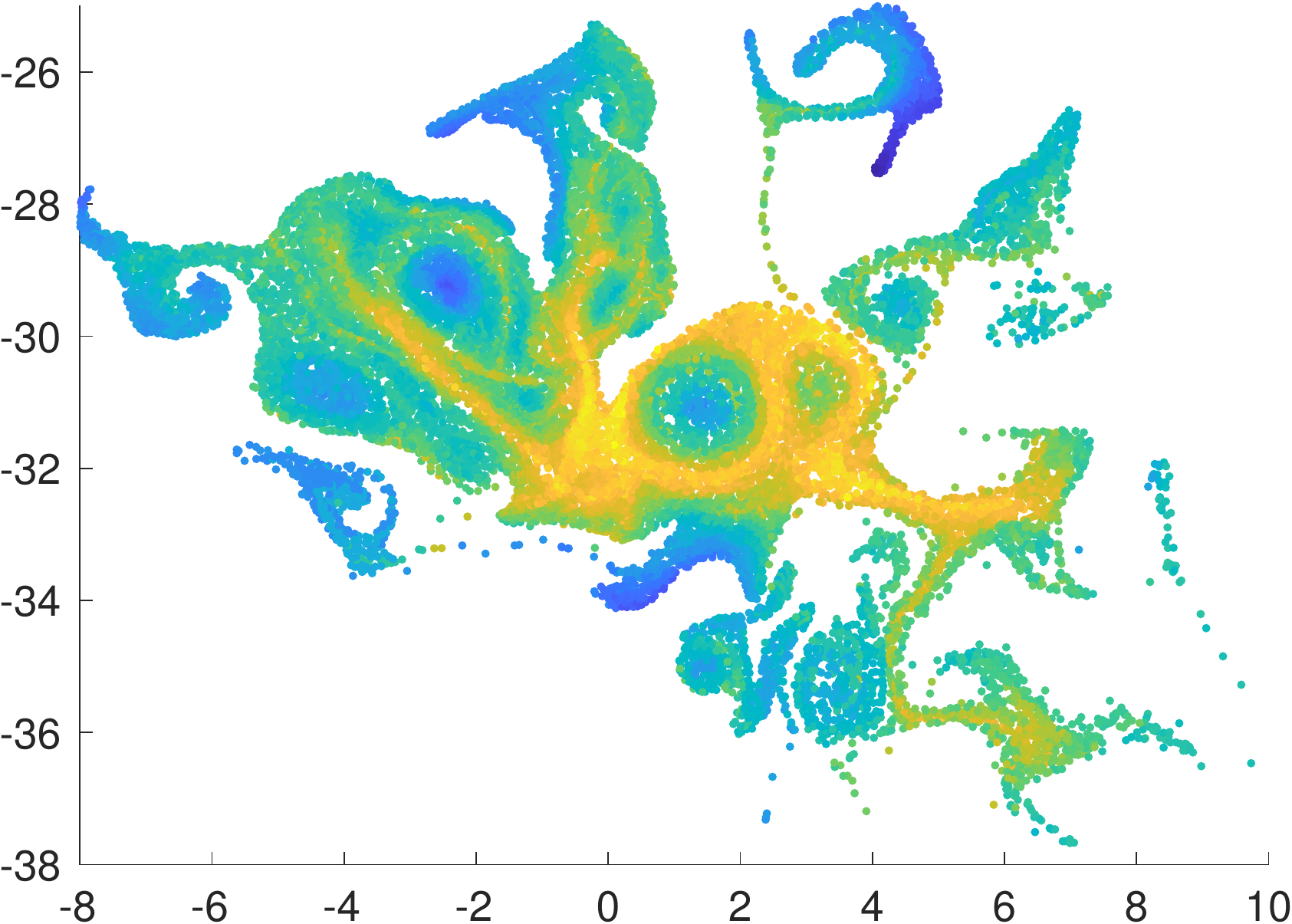}
\hfill
\includegraphics[width = 0.24\textwidth]{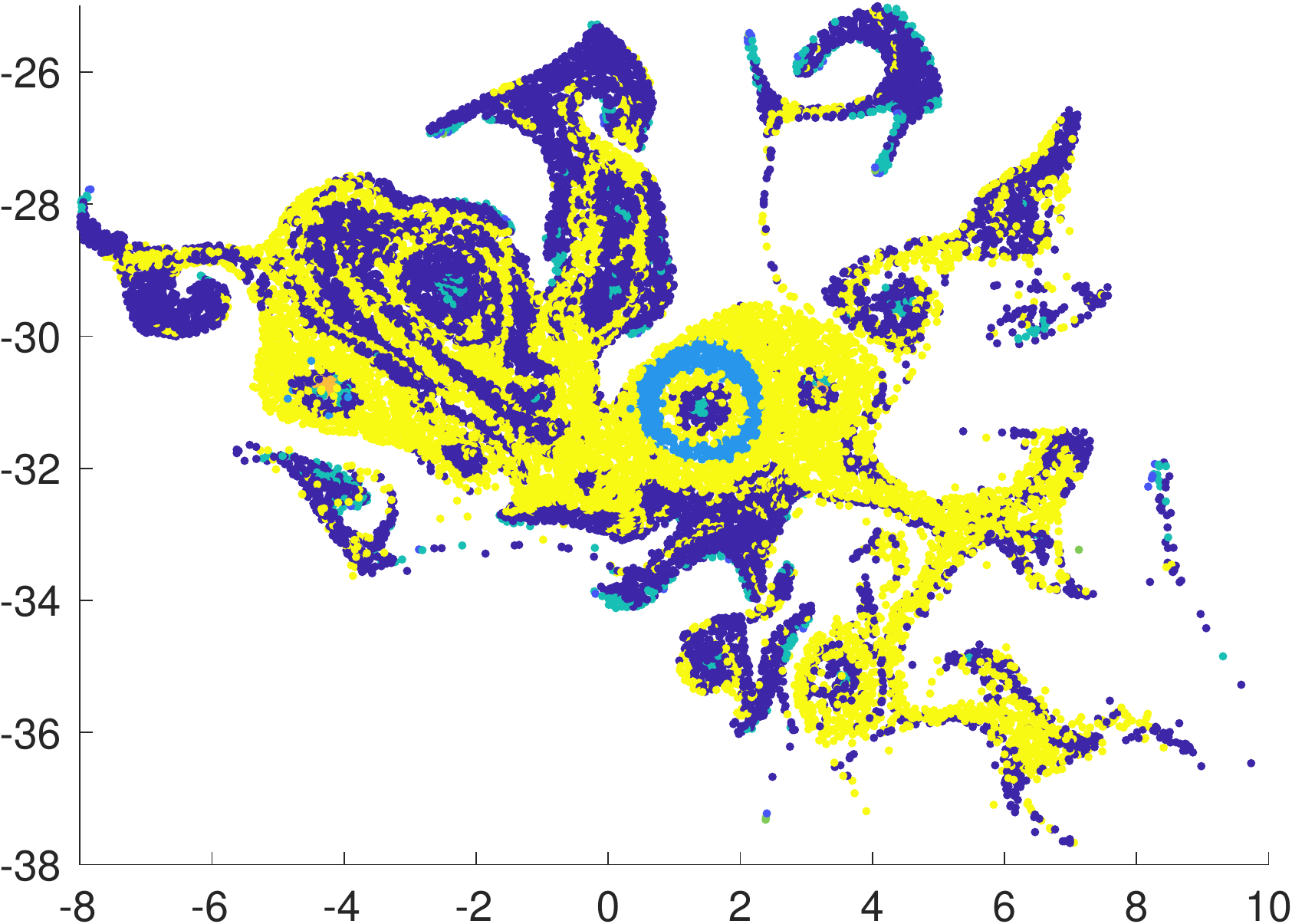}
\caption{Ocean flow data set. From left to right: degree, clustering coefficient, closeness, classification. Top: initial time, bottom: final time. For the first three columns lighter yellow colors indicate higher values. while darker blue indicate lower ones. In the rightmost column the colors indicate seven clusters.}
\label{fig:aviso1}
\end{figure*}

As a last experimental case, we will now analyze an actual ocean flow. We consider a velocity field of the surface water derived from AVISO satellite altimetry measurements. The flow is area-preserving on a spherical surface. We focus on the region of the Agulhas leakage in the South Atlantic Ocean, using the same data set as refs.\ \onlinecite{HaHa2016, FrJu18}. We initialize a $200\times 100$ array of drifters in $[-4,6]\times [-34,-28]$ advected by the flow and observed at times~$t=0,1,\ldots,90$, where $t_0$ corresponds to November 11, 2006 \cite{HaHa2016,FrJu18}. With $\ep=0.1$ we calculate the adjacency matrix and the so far discussed network measures, then compute a classification of the trajectories based on degree and clustering coefficient, as done in the previous example (i.e., we perform a clustering of the diffusion-map embedding of the degree-clustering coefficient point cloud). The results are shown and described  in Figure~\ref{fig:aviso1}.

First, we observe that the highest degree attained---in contrast to the double gyre flow---is on the outer perimeter of an eddy. This is underlined by closeness, as the neighborhood of this eddy seems to be some sort of ``hub'' for transport; many trajectories from different regions pass by this eddy. Second, we also observe trajectories of high clustering coefficient value near the (time-evolving) boundary of the region of consideration. This is due to the dynamics-induced filamentation; subsets of trajectories are separated from the ``main region'' and build islands or peninsula that do not return to an~$\ep$-proximity of other trajectories. Thus, this subset maintains a low degree and high internal connectivity, giving a large clustering coefficient value.

Both of these behaviors arise, because this flow, in contrast to the previous ones, is considered on a ``free domain''; meaning that the state space (the region we have trajectory data from) evolves with the flow. Thus, we are not taking dynamical information into account from the neighborhood of our set of trajectories. However this neighborhood interacts with our observations, as the region where we have trajectory information starts to mix with the white region, where we do not have any. In summary, this ``free domain'' situation pollutes our analysis with spurious structures. The derivation of sensible network measures that account for this dynamical situation is the next challenge on the way to being able to apply these methods in real world situations.

\section{Conclusion}
\label{sec:conclusion}
We have studied an unweighted and undirected trajectory-based network\cite{PGSc17}. Simply computable network measures allow us to infer valuable information about the dynamics of the underlying system \cite{faranda2018correlation}, even if a full global analysis of the system is out of reach---due to the dimensionality of the system, or because only a finite amount of trajectory data is available.

The palette of network measures is broad, both in complexity and computational efficiency. We have focused on simpler ones here, and were able to show analytic connections between the local degree of a network  and quantitative dynamical descriptors, like FTLE, in the large-data limit.

More complex network measures, such as clustering coefficient and closeness can be linked to qualitative dynamical behavior. We have experimentally verified these connections, and have shown how classification with respect to multiple network measures separates regions exhibiting different (topological) dynamical behavior.

The general aim is to identify structurally different dynamical behavior from large sets of possibly high-dimensional trajectory data. Further developments need to be done to understand how consistent estimators of dynamical descriptors can be derived from trajectory networks, how to deal with the ``free domain problem'' above, and with missing data, in general.

\section*{Acknowledgments}
This work is supported by the Deutsche Forschungsgemeinschaft (DFG) through the Priority Programme SPP 1881 ``Turbulent Superstructures''.
PK also acknowledges funding from the Deutsche Forschungsgemeinschaft (DFG) through the Collaborative Research Center ``Scaling Cascades of Complex Systems'', project A01.
KPG also acknowledges funding from EU Marie-Sk{\l}odowska-Curie ITN Critical Transitions in Complex Systems (H2020-MSCA-2014-ITN 643073 CRITICS). 

\vspace*{5mm}
\bibliographystyle{abbrv}
\bibliography{library}

\end{document}